\documentclass[paper,twocolumn,twoside]{geophysics}
\usepackage{amssymb}
\usepackage{amsmath}
\usepackage{amscd}
\usepackage{graphicx}
\usepackage{graphics}
\usepackage[noadjust]{cite}
\usepackage{amsthm}
\usepackage{caption}
\usepackage{multirow}
\usepackage{array}
\usepackage{rotating}
\usepackage[norelsize,ruled,linesnumbered]{algorithm2e}
\usepackage{algorithm2e,setspace}

\usepackage{indentfirst}
\usepackage{tikz}
\usepackage{calc}
\usepackage{epsfig}
\usepackage{natbib}
\usepackage[makeroom]{cancel}
\usepackage{multicol}
\usepackage{csquotes}

\newcommand{\ket}{\right\rangle}       
\newcommand{\bra}{\left\langle}

\newsavebox\CBox
\newcommand\hcancel[2][0.5pt]{%
  \ifmmode\sbox\CBox{$#2$}\else\sbox\CBox{#2}\fi%
  \makebox[0pt][l]{\usebox\CBox}%
  \rule[0.5\ht\CBox-#1/2]{\wd\CBox}{#1}}

\newcommand\Tstrut{\rule{0pt}{20pt}}         
\begin{document}
\title{ADMM-based multi-parameter wavefield reconstruction inversion in VTI acoustic media with TV regularization}

\address{ \footnotemark[1]University of Tehran, Institute of Geophysics, Tehran, Iran, email: h.aghamiry@ut.ac.ir, agholami@ut.ac.ir \\ \\
\footnotemark[2]University Cote d'Azur - CNRS - IRD - OCA, Geoazur, Valbonne, France, email: aghamiry@geoazur.unice.fr, operto@geoazur.unice.fr}
\author{Hossein S. Aghamiry \footnotemark[1]\footnotemark[2], Ali Gholami \footnotemark[1] and St\'ephane Operto \footnotemark[2]}

\lefthead{A PREPRINT~~~-~~~Aghamiry et al.}
	\righthead{~~~~~~~~~~~~~~~~~~~~~~~~~~~~A PREPRINT~~~-~~~IR-WRI for VTI acoustic medium}
\maketitle
\begin{abstract}
\footnotesize{
Full waveform inversion (FWI) is a nonlinear waveform matching procedure, which suffers from cycle skipping when the initial model is not kinematically-accurate enough. To mitigate cycle skipping, wavefield reconstruction inversion (WRI) extends the inversion search space by computing wavefields with a  relaxation of the wave equation in order to fit the data from the first iteration. Then, the subsurface parameters are updated by  minimizing the source residuals the relaxation generated.
Capitalizing on the wave-equation bilinearity, performing wavefield reconstruction and parameter estimation in alternating mode decomposes WRI into two linear subproblems, which can solved efficiently with the alternating-direction method of multiplier (ADMM), leading to the so-called IR-WRI. Moreover, ADMM provides a suitable framework to implement bound constraints and different types of regularizations and their mixture in IR-WRI. 
Here, IR-WRI is extended to multiparameter reconstruction for VTI acoustic media. To achieve this goal, we first propose different forms of bilinear VTI acoustic wave equation. We develop more specifically IR-WRI for the one that relies on a parametrisation involving vertical wavespeed and Thomsen's parameters $\delta$ and $\epsilon$. With a toy numerical example, we first show that the radiation patterns of the virtual sources generate similar wavenumber filtering and parameter cross-talks in classical FWI and IR-WRI. Bound constraints and TV regularization in IR-WRI fully remove these undesired effects for an idealized piecewise constant target. We show with a more realistic long-offset case study representative of the North Sea that anisotropic IR-WRI successfully reconstruct the vertical wavespeed starting from a laterally homogeneous model and update the long-wavelengths of the starting $\epsilon$ model, while a smooth $\delta$ model is used as a passive background model. VTI acoustic IR-WRI can be alternatively performed with subsurface parametrisations involving stiffness or compliance coefficients or normal moveout velocities and $\eta$ parameter (or horizontal velocity)}.
\end{abstract}
%
\section{Introduction}
Full waveform inversion (FWI) is a waveform matching procedure which provides subsurface model with a wavelength resolution. However, it suffers from cycle skipping when the initial model is not  accurate enough according to the lowest available frequency.
To mitigate cycle skipping, the search space of frequency-domain FWI can be extended by wavefield reconstruction inversion (WRI) \citep{VanLeeuwen_2013_MLM,vanLeeuwen_2016_PMP}. In WRI, the search space is extended with a penalty method to relax the wave-equation constraint at the benefit of the observation equation (i.e., the data fit) during wavefield reconstruction. Then, the subsurface parameters are estimated from the reconstructed wavefields by minimizing the source residuals the relaxation generated. If these two subproblems (wavefield reconstruction and parameter estimation) are solved in alternating mode \citep{VanLeeuwen_2013_MLM} rather than by variable projection \citep{vanLeeuwen_2016_PMP}, WRI can be recast as a sequence of two linear subproblems capitalizing on the bilinearity of the scalar Helmholtz equation with respect to the wavefield and the squared slownesses (the Helmholtz equation is linear with respect to the wavefield for a given squared slowness model and is linear with respect to the squared slownesses for a given wavefield).  \\
\citet{Aghamiry_2019_IWR} improved WRI by replacing the penalty method with the augmented Lagrangian method \citep{Nocedal_2006_NOO} and solve it using the alternating-direction method of multipliers (ADMM) \citep{Boyd_2011_DOS}, leading to iteratively-refined WRI (IR-WRI). Although ADMM was originally developed for separable convex problems, the bilinearity of the wave equation makes IR-WRI biconvex, which allows for the use of ADMM as is \citep[][ Section 3.1.1]{Boyd_2011_DOS}. Moreover, a scaled form of ADMM draws clear connection between WRI and IR-WRI in the sense that it shows that IR-WRI reduces to a penalty method in which the right-hand sides (RHSs) in the quadratic objective functions associated with the observation equation and the wave equation are iteratively updated with the running sum of the data and source residuals in iteration (namely, the scaled Lagrange multipliers). This RHS updating, the lacking feature in WRI, efficiently refines the solution of the two linear subproblems when a fixed penalty parameter is used \citep{Aghamiry_2019_IWR}.  Later, \citet{Aghamiry_2018_MIA,Aghamiry_2019_IBC} interfaced bound constraints and edge-preserving regularizations with ADMM to manage large-contrast media. Also, to preserve the smooth components of the subsurface when edge preserving regularizations are used, \citet{Aghamiry_2018_HTT,Aghamiry_2019_CROb} combine blocky and smooth promoting regularization in the framework of IR-WRI by using infimal convolution of Tikhonov and Total Variation (TV) regularization functions. 

IR-WRI was mainly assessed for wavespeed estimation from the scalar Helmholtz equation. The objective of this study is to develop and assess the extension of IR-WRI to multi-parameter reconstruction in VTI acoustic media. 
To achieve this goal, we first need to review different formulations of the VTI acoustic wave equation and different subsurface parametrisations for which bilinearity of the wave equation is fulfilled, in order to keep the parameter estimation subproblem linear. 
Since the wavefield reconstruction requires to solve a large and sparse system of linear equations, second-order or fourth-order wave equation for pressure will be favored at the expense of first-order velocity-stress formulations to mitigate the number of unknowns during wavefield reconstruction. However, we stress that different forms of the wave equation can be used to perform wavefield reconstruction and parameter estimation, provided they provide consistent solutions \citep[][ Their Appendix A and B]{Gholami_2013_WPA1}. \\
Also, we will favor subsurface parametrisation involving the vertical wavespeed $\bold{v}_0$ and the Thomsen's parameter $\epsilon$ and $\delta$ at the expense of that involving stiffness coefficients according to the trade-off analysis of \citet[][ Their Appendix A and B]{Gholami_2013_WPA1}. 

When bilinearity of the wave equation is  fulfilled, IR-WRI can be extended to multi-parameter estimation following the procedure promoted by \citet{Aghamiry_2019_IBC}, where TV regularization and bound constraints are efficiently implemented in the parameter-estimation subproblem using variable splitting schemes \citep{Glowinski_2017_SMC}. The splitting procedure allows us to break down the non-differentiable TV regularization problem into two easy-to-solve subproblems: a least-squares quadratic subproblem and a proximity subproblem \citep{Goldstein_2009_SBM}. 

In this study, we perform a first assessment of multi-parameter IR-WRI for $\bold{v}_0$ and $\epsilon$ using a toy inclusion example and the more realistic synthetic North Sea case study tackled by \citet{Gholami_2013_WPA1} and \citet{Gholami_2013_WPA2}.
With the toy example, we first show that the radiation patterns and the parameter cross-talks have the same effects as in classical FWI when TV regularization is not applied. Then, we show how TV regularization fully removes bandpass filtering and cross-talk effects generated by these radiation patterns for this idealized piecewise constant target.
With the North Sea  example, we first show the resilience of IR-WRI against cycle skipping when using a crude initial $\bold{v}_0$ model. The reconstruction of the $\bold{v}_0$ model is accurate except in the deep smooth part of the subsurface which suffers from a deficit of wide-angle illumination, while the reconstruction of $\epsilon$ is more challenging and requires additional regularization to keep its update smooth and close to the starting model. However, we manage to update significantly the long-to-intermediate wavelengths of $\epsilon$ by IR-WRI, unlike \citet{Gholami_2013_WPA2}. Also, comparison between mono-parameter IR-WRI for $\bold{v}_0$ and multi-parameter IR-WRI for $\bold{v}_0$ and $\epsilon$ allows us to gain qualitative insights on the sensitivity of the IR-WRI to $\epsilon$ in terms of subsurface model quality and data fit.

This paper is organized as follow. We first discuss the bilinearity of the acoustic VTI wave equation as well as its implication on the gradient and the Gauss-Newton Hessian of the parameter-estimation subproblem. From the selected bilinear formulation of the wave equation, we develop bound-constrained TV-regularized IR-WRI for VTI acoustic media parametrized by $1/\bold{v}_0^2$, $1 + 2 \epsilon$ and $\sqrt{1 + 2 \delta}$. Third, we assess IR-WRI for VTI acoustic media against the inclusion and North Sea  case studies. Finally, we discuss the perspectives of this work.


\section{Theory}
In this section, we first show that the VTI acoustic wave equations is bilinear with respect to the wavefield and model parameters. 
Then, we rely on this bilinearity to formulate multi-parameter IR-WRI in VTI acoustic media with bounding constraints and TV regularization.

\subsection{Bilinearity of the wave equation: preliminaries}
If we write the wave equation in a generic matrix form as
\begin{equation}
\bold{A}(\bold{m}) \bold{u} = \bold{s},
\label{eqAus}
\end{equation}
where $\bold{A} \in \mathbb{C}^{(n \times n_c) \times (n \times n_c)}$ is the wave-equation operator, $\bold{u} \in \mathbb{C}^{(n \times n_c) \times 1}$ is the wavefield vector, $\bold{s} \in \mathbb{C}^{(n \times n_c) \times 1}$ is the source vector, $\bold{m} \in \mathbb{R}^{(n \times n_m) \times 1}$ is the subsurface parameter vector, $n$ is the number of degrees of freedom in the spatial computational mesh, $n_c$ is the number of wavefield components, $n_m$ is the number of parameter classes, then
the wave equation is bilinear if there exists a linear operator $\bold{L}(\bold{u}) \in \mathbb{C}^{(n \times n_m) \times (n \times n_m)}$ such that
\begin{equation}
\bold{L}(\bold{u}) \bold{m} = \bold{y}(\bold{u}),
\label{eqLmy}
\end{equation}
where $\bold{y}(\bold{u}) \in \mathbb{C}^{(n \times n_m) \times 1}$.
%
Bilinearity is verified when the left-hand side of the wave equation can be decomposed as:
\begin{equation}
\bold{A}(\bold{m}) \bold{u} = \bold{B} \, \bold{M}(\bold{m}) \, \bold{C} \bold{u} + \bold{D} \bold{u},
\label{eqwe}
\end{equation}
where $\bold{M}(\bold{m})$ is a block matrix, whose blocks of dimension $n \times n$ are either 0 or of the form $\text{diag}(\bold{m}_i)$, and  matrices $\bold{B}$, $\bold{C}$ and $\bold{D}$ don't depend on $\bold{m}$. The operator $\text{diag}(\bullet)$ denotes a diagonal matrix of coefficients $\bullet$ and $\bold{m}_i$ is the subsurface parameter vector of class $i$. 
By noting that $\text{diag}(\bold{x}) \bold{z} = \text{diag}(\bold{z}) \bold{x}$, the block diagonal structure of $\bold{M}$ allows one to rewrite the term $\bold{B} \, \bold{M}(\bold{m}) \, \bold{C} \bold{u} $ as $\bold{B} \, \bold{C}'(\bold{u}) \bold{m}$, where $\bold{C}'(\bold{u})$ is a block matrix, whose blocks of dimension $n \times n$ are either 0 or diagonal with coefficients depending on $\bold{u}$.

It follows from equation \ref{eqwe} and the above permutation between $\bold{u}$ and $\bold{m}$ that the wave equation can be re-written as
\begin{equation}
\bold{A}(\bold{m}) \bold{u} - \bold{s} = \bold{B} \, \bold{C}'(\bold{u}) \bold{m} + \bold{D} \bold{u} - \bold{s} = \bold{L}(\bold{u}) \bold{m} - \bold{y}(\bold{u}).
\end{equation}
Moreover, in the framework of multi-parameter analysis, it is worth noting that 
\begin{equation}
\frac{\partial \bold{A}(\bold{m})}{\partial m_{k}} \bold{u} = \bold{L}(\bold{u}) \bold{e}_{k},
\end{equation}
where the left-hand side is the so-called virtual source associated with $m_k$ \citep{Pratt_1998_GNF} and $\bold{e}_{k} \in \mathbb{C}^{(n \times n_m) \times 1}$ denotes a column vector whose $k$th component is one while all the others are zeros.

Accordingly, the normal operator $\bold{L}^T \bold{L}$, i.e., the Gauss-Newton Hessian of the parameter estimation subproblem in IR-WRI, is formed by the zero-lag correlation of the virtual sources
%
%
and, hence is extremely sparse. We also point that $\bold{L}$ and its associated normal operator $\bold{L}^T \bold{L}$ are block diagonals if $\bold{B}$ is diagonal. This means that, if the model parameters are first sorted according to their position in the mesh (fast index) and second according to the parameter class they belong to (slow index), then the diagonal coefficients of the off-diagonal blocks describe the inter-parameter coupling.
%
%
In the following section, we show the bilinearity of the acoustic VTI wave equation based upon the above matrix manipulations.
%
%
\subsection{Bilinearity of the acoustic VTI wave equation}
\subsection{First-order velocity-stress wave equation}
We first consider the frequency-domain first-order velocity-stress wave equation in 2D VTI acoustic media \citep{Duveneck_2008_AVW,Duveneck_2011_SPM,Operto_2014_FAT} 
\begin{align} \label{edavti}
-\hat{i} \omega \bold{v}_{x,l} &= \text{diag}(\bold{b})  \nabla_{\!\! x} \bold{u}_{x,l},  \nonumber \\
-\hat{i} \omega \bold{v}_{z,l} &= \text{diag}(\bold{b}) \nabla_{\!\! z}\bold{u}_{z,l},   \nonumber \\
-\hat{i} \omega \bold{u}_{x,l} &= \text{diag}(\bold{c}_{11}) \nabla_{\!\! x}\bold{v}_{x,l}  + \text{diag}(\bold{c}_{13}) \nabla_{\!\! z} \bold{v}_{z,l} - \hat{i} \omega \bold{s}_l, \nonumber \\
-\hat{i} \omega \bold{u}_{z,l} &= \text{diag}(\bold{c}_{13}) \nabla_x \bold{v}_{x,l}  + \text{diag}(\bold{c}_{33}) \nabla_{\!\! z} \bold{v}_{z,l} - \hat{i} \omega \bold{s}_l, 
\end{align}
where $\hat{i}=\sqrt{-1}$, $\omega$ is the angular frequency, $\bold{v}_{x,l} \in \mathbb{C}^{n \times 1}$ and $\bold{v}_{z,l} \in \mathbb{C}^{n \times 1}$ are the horizontal and vertical particle velocity wavefields, $\bold{s}_l \in \mathbb{C}^{n \times 1}$ denote the pressure sources, and $\bold{u}_{x,l}\in \mathbb{C}^{n \times 1}$ and $\bold{u}_{z,l}\in \mathbb{C}^{n \times 1}$ are the so-called horizontal and vertical pressure wavefields \citep{Plessix_2011_GJI}. The subscript $l \in \{1,2,...,n_s\}$ is the source index, where $n_s$ denotes the number of sources.
The subsurface properties are parametrized by the buoyancy $\bold{b} \in \mathbb{R}^{n \times 1}$ (inverse of density) and the stiffness coefficients $\bold{c}_{ij} \in \mathbb{R}^{n \times 1}$.
Operators $\nabla_{\!\! x}$ and $\nabla_{\!\! z}$ are finite difference approximation of first order derivative operators with absorbing perfectly matched layer (PML) coefficients \citep{Berenger_1994_PML}.

Gathering equation \ref{edavti} for all sources leads to the following matrix equation:
\begin{equation}
\label{eqvelstressf5} 
 \left( 
 \begin{bmatrix}
\bold{M}_b & \bold{0}\\
\bold{0} & \bold{M}_c\\
\end{bmatrix} 
 \begin{bmatrix}
\bold{0} & \bold{\nabla}\\
\bold{\nabla} & \bold{0}\\
\end{bmatrix} 
  + \hat{i} \omega  \bold{I}\right)
 \begin{bmatrix}
\bold{V}\\
\bold{U}\\
\end{bmatrix} 
 =   \hat{i}\omega
  \begin{bmatrix}
\bold{0}\\
\bold{S}\\
\end{bmatrix}, 
\vspace{-0.3cm}
\end{equation}
%
%
where $\bold{I}$ is the identity matrix,
\begin{equation}
\nabla=
\begin{bmatrix}
\nabla_{\!\! x}  & \bold{0}\\
\bold{0}  & \nabla_{\!\! z} \\
\end{bmatrix}, \nonumber
\end{equation}
\begin{equation} \label{UVS}
\bold{U}=
 \begin{bmatrix}
   \bold{U}_h  \\
   \bold{U}_v
\end{bmatrix}, \quad
\bold{V}=
 \begin{bmatrix}
   \bold{V}_x  \\
   \bold{V}_z
\end{bmatrix}, \quad
\bold{S}=
\begin{bmatrix}
      \bold{S}_x  \\
      \bold{S}_z
\end{bmatrix}, 
\end{equation} 
%
with 
\begin{equation}
\bold{V}_x= 
\begin{bmatrix}
\bold{v}_{x,1} & \bold{v}_{x,2} & ... & \bold{v}_{x,n_s}\
\end{bmatrix} \in \mathbb{C}^{n \times n_s}, \nonumber
\end{equation}
and analogously for $\bold{V}_z$, $\bold{U}_h$, $\bold{U}_v$, and
\begin{equation}
\bold{S}_{x} = \bold{S}_{z}=
\begin{bmatrix}
\bold{s}_{1} & \bold{s}_{2} & ... & \bold{s}_{n_s}
\end{bmatrix} \in \mathbb{C}^{n \times n_s}. \nonumber
\end{equation}
Furthermore,
\begin{equation*}
\bold{M}_b=
\begin{bmatrix}
\text{diag}(\bold{b})  & \bold{0} \\
\bold{0}  & \text{diag}(\bold{b}) 
\end{bmatrix},
\bold{M}_c=
\begin{bmatrix}
\text{diag}(\bold{c}_{11}) & \text{diag}(\bold{c}_{13}) \\
\text{diag}(\bold{c}_{13}) & \text{diag}(\bold{c}_{33}) \\
\end{bmatrix}. \nonumber
\end{equation*}

%

Note that, according to the decomposition introduced in the previous section, 
$\bold{M}= \begin{bmatrix}
\bold{M}_b & \bold{0}  \\
\bold{0} & \bold{M}_c
\end{bmatrix}$, $\bold{B}=\bold{I}$,  $\bold{C}=\nabla$ and $\bold{D}=\hat{i} \omega \bold{I}$.
Equation \ref{eqvelstressf5} is linear in $\bold{U}$ and $\bold{V}$ when the model parameters embedded in $\bold{M}_b$ and $\bold{M}_c$  are known. 
When $\bold{U}$ and $\bold{V}$  are known, this system can be also recast as a new linear system in which the unknowns are the model parameters, thus highlighting the bilinearity of the wave equation.
For the $l$th source, the new equations become
\begin{equation} \label{Eqval_m}
\begin{bmatrix}
\bold{L}_{11} & \bold{0} & \bold{0} & \bold{0}   \\
\bold{L}_{21} & \bold{0} & \bold{0} & \bold{0} \\
 \bold{0} & \bold{L}_{32} & \bold{L}_{33}& \bold{0}  \\
 \bold{0} & \bold{0} & \bold{L}_{43}  & \bold{L}_{44}
\end{bmatrix}
\begin{bmatrix}
\bold{b} \\
\bold{c}_{11} \\
\bold{c}_{13} \\
\bold{c}_{33} \\
\end{bmatrix}
= \hat{i} \omega 
\begin{bmatrix}
\bold{0}-\bold{v}_{x,l} \\
\bold{0}-\bold{v}_{z,l} \\
 {\bold{s}_l}-\bold{u}_{x,l}\\
 {\bold{s}_l}-\bold{u}_{z,l}\\
\end{bmatrix},
\end{equation}
where
\begin{equation}
\begin{cases}
\bold{L}_{11}=\text{diag}(\nabla_{\!\! x} \bold{u}_{x,l}), \\
 \bold{L}_{21}=\text{diag}(\nabla_{\!\! z} \bold{u}_{z,l}), \\
\bold{L}_{32}=\bold{L}_{43}=\text{diag}(\nabla_{\!\! x} \bold{v}_{x,l}),\\
\bold{L}_{33}=\bold{L}_{44}=\text{diag}(\nabla_{\!\! z} \bold{v}_{z,l}).
 \end{cases} \nonumber
\end{equation}
Equations \ref{eqvelstressf5} and \ref{Eqval_m} are equivalent forms of the original equation \ref{edavti}.
The former expresses the discretized wavefields as the unknowns of a linear system, whose coefficients depend on the known subsurface parameters, while the latter expresses the model parameters as the unknowns of an another linear system, whose coefficients depend on the known wavefields.
In the framework of WRI, this bilinearity allows one to recast the waveform inversion problem as two linear subproblems for wavefield reconstruction and parameter estimation, which can be solved efficiently in alternating mode with ADMM \citep{Aghamiry_2019_IWR,Aghamiry_2019_IBC}.
In the next section, we show the bilinearity of the second-order frequency-domain wave equation, which may be more convenient to solve with linear algebra methods than the first-order counterpart, since it involves fewer unknowns for a computational domain of given size.
\subsection{Second-order wave equation}
Following a parsimonious approach \citep[e.g.,][]{Operto_2007_3D}, we eliminate $\bold{v}_{x,l}$ and $\bold{v}_{z,l}$ from equation \ref{edavti} to derive a system of two second-order partial differential equations as
\begin{equation} \label{sovti}
\begin{split}
& - \omega^2 \bold{u}_{x,l} = \text{diag}(\bold{c}_{11}) \nabla_{\!\! x} \text{diag}(\bold{b})  \nabla_{\!\! x}  \bold{u}_{x,l}  \\
& \hspace{1.4cm}+ \text{diag}(\bold{c}_{13}) \nabla_{\!\! z}\text{diag}(\bold{b})  \nabla_{\!\! z} \bold{u}_{z,l} -  \omega^2  \bold{s}_l, \\
& - \omega^2 \bold{u}_{z,l} = \text{diag}(\bold{c}_{13}) \nabla_{\!\! x} \text{diag}(\bold{b}) \nabla_{\!\! x} \bold{u}_{x,l} \\
& \hspace{1.4cm} + \text{diag}(\bold{c}_{33}) \nabla_{\!\! z} \text{diag}(\bold{b})  \nabla_{\!\! z} \bold{u}_{z,l} - \omega^2  \bold{s}_l.
\end{split}
\end{equation}
Equation~\ref{sovti} defines a tri-linear equation with respect to buoyancy, stiffness parameters and pressure wavefields. A tri-linear function is a function of three variables which is linear in one variable when the other two variables are fixed. 
In this study, we will assume that density is constant and equal to 1 to focus on the estimation of the anisotropic parameters. If heterogeneous density needs to be considered, the second-order wave equation can be recast as a bilinear system if the first-order wave equation is parametrized with compliance coefficients instead of stiffness coefficients \citep[see Appendix A and ][]{Vigh_2014_EFI,Yang_2016_SFM}. Alternatively, the parameter estimation can be performed with the bilinear first-order wave equation, equation \ref{Eqval_m}, while wavefield reconstruction is performed with the second-order wave equation, equation \ref{sovti}, or the fourth-order wave equation reviewed in Appendix C for sake of computational efficiency. 

We continue by assuming that the density is constant and equal to 1 and parametrize the VTI equation in terms of vertical wavespeed $\bold{v}_0$ and Thomsen's parameters $\bold{\epsilon}$ and $\bold{\delta}$ \citep{Thomsen_1986_WEA}.
Accordingly, we rewrite equation \ref{sovti} as
\begin{equation} \label{sovtit}
\begin{split}
& - \omega^2 \text{diag}\left(1/\bold{v}_{0}^2\right) \bold{u}_{x,l} = \text{diag}(1 + 2 \epsilon) \nabla_{\!\! xx}  \bold{u}_{x,l}  \\
& \hspace{1.4cm}+ \text{diag}(\sqrt{1 + 2 \delta})  \nabla_{\!\! zz} \bold{u}_{z,l} - \text{diag}(1/\bold{v}_{0}^2)  \bold{s}_l, \\
& -  \omega^2 \text{diag}\left(1/\bold{v}_{0}^2\right) \bold{u}_{z,l} = \text{diag}(\sqrt{1 + 2 \delta}) \nabla_{\!\! xx}  \bold{u}_{x,l} \\
& \hspace{1.4cm} + \nabla_{\!\! zz}  \bold{u}_{z,l} - \text{diag}\left(1/\bold{v}_{0}^2\right)  \bold{s}_l,
\end{split}
\end{equation}
where $\nabla_{\!\! xx}=\nabla_{\!\! x}\nabla_{\!\! x}$ and $\nabla_{\!\! zz}=\nabla_{\!\! z}\nabla_{\!\! z}$. 
We write this $2n \times 2n$ linear system in a more compact form as
\begin{equation}
\label{linearu}
\bold{A}(\bold{m}) \bold{U} = \bold{S},
\end{equation}
where 
$\bold{S}^T= \omega^2
\begin{bmatrix}
      \bold{S}^T_x  &
      \bold{S}^T_z
\end{bmatrix}$,
$\bold{U}$ is defined as in equation \ref{UVS}, and
the matrix $\bold{A}$ is given by
\begin{equation*}
\begin{bmatrix}
    \omega^2  \text{diag}(\bold{m}_{v_0})+ \text{diag}(\bold{m_\epsilon}) \nabla_{\!\! xx} & \text{diag}(\bold{m_\delta}) \nabla_{\!\! zz} \\
    \text{diag}(\bold{m_\delta}) \nabla_{\!\! xx} & \omega^2  \text{diag}(\bold{m}_{v_0}) + \nabla_{\!\! zz} 
\end{bmatrix},
\end{equation*}
where the model parameters are
\begin{equation}
\bold{m}=
 \begin{bmatrix}
   \bold{m}_{v_0}  \\
   \bold{m_{\epsilon}}   \\
   \bold{m_{\delta}} 
\end{bmatrix}=
 \begin{bmatrix}
   1/{\bold{v}_{0}^2}  \\
  1+2\epsilon   \\
  \sqrt{1+2\bold{\delta}}
\end{bmatrix} \nonumber
.
\end{equation}
%
%
Equation \ref{linearu} is linear in wavefields $\bold{U}$ when the model parameters $\bold{m}$ are known. 
When $\bold{U}$ is known, this system can be recast as a linear system in which the unknowns are the model parameters.
\begin{equation}
\label{linearma}
\begin{bmatrix}
\bold{L}_{1} \\
\vdots \\
\bold{L}_{l} \\
\vdots \\
\bold{L}_{n_s}
\end{bmatrix}
\bold{m} = 
\begin{bmatrix}
\bold{y}_1 \\
\vdots \\
\bold{y}_l \\
\vdots \\
\bold{y}_{n_s}
\end{bmatrix},
\end{equation}
%
where $\bold{L}_l$ is given by
\begin{equation*}
\begin{bmatrix}
    \omega^2  \text{diag}(\bold{u}_{x,l}) & \text{diag}(\nabla_{\!\! xx} \bold{u}_{x,l})  & \text{diag}(\nabla_{\!\! zz} \bold{u}_{z,l}) \\
   \omega^2   \text{diag}(\bold{u}_{z,l}) & \bold{0} &  \text{diag}(\nabla_{\!\! xx} \bold{u}_{x,l})    
\end{bmatrix}
\end{equation*}
and
\begin{equation}
\bold{y}_{l}=
 \begin{bmatrix}
    \omega^2 \bold{s}_l  \\
    \omega^2 \bold{s}_l -\nabla_{\!\! zz} \bold{u}_{z,l}
\end{bmatrix}. \nonumber
\end{equation}
Note that each block of $\bold{L}_l$ is diagonal.
In the next section, we develop multi-parameter acoustic VTI IR-WRI with bound constraints and TV regularization. 
We give the most general formulation in which all the tree parameter classes $\bold{m}_{v_0}$, $\bold{m_{\epsilon}}$, $\bold{m_{\delta}}$ 
are optimization parameters (updated by the inversion). However, one may process some of them as passive parameters or update the parameter classes in sequence rather than jointly.
In this case, the linear system associated with the parameter estimation subproblem, equation \ref{linearma}, will change. 
Table \ref{Tab1} presents this system for the different possible configurations.
\begin{table*}[ht!]
\captionof{table}{The linear systems corresponding to the update of $\bold{m}_{v_0}$, $\bold{m}_\epsilon$ and $\bold{m}_\delta$ when they are active or passive during the inversion. 
In the first three columns, \enquote{a.} denotes an active parameter and \enquote{p.} denotes a passive parameter.}
\label{Tab1}
\begin{center}
\begin{tabular}{|c|c|c|c|}
 \hline 
$\bold{m}_{v_0}$ & $\bold{m}_\epsilon$ & $\bold{m}_\delta$ & The corresponding linear system for updating the model \Tstrut \\ [5pt] \hline
 a. & p. & p. & $\underbrace{
\begin{bmatrix}
   \omega^2  \text{diag}(\bold{u}_{x,l}) \\
   \omega^2  \text{diag}(\bold{u}_{z,l})      
\end{bmatrix}}_{\bold{L}_{l} \in \mathbb{C}^{ 2n  \times  n }}
\underbrace{ \begin{bmatrix}
   \bold{m}_{v_0}   
\end{bmatrix}}_{\bold{m} \in \mathbb{R}^{n \times 1}}
=
 \underbrace{
 \begin{bmatrix}
    \omega^2 \bold{s}_l- \text{diag}(\bold{m}_\epsilon) \nabla_{\!\! xx} \bold{u}_{x,l}  - \text{diag}(\bold{m}_\delta)\nabla_{\!\! zz} \bold{u}_{z,l}\\
    \omega^2 \bold{s}_l -\nabla_{\!\! zz} \bold{u}_{z,l} - \text{diag}(\bold{m}_\delta){\nabla_{\!\! xx} \bold{u}_{x,l}}
\end{bmatrix}
}_{\bold{y}_{l} \in \mathbb{C}^{2n \times 1}}$ \Tstrut \\ [25pt] \hline
 p. & a. & p. & $\underbrace{
\begin{bmatrix}
    \text{diag}(\nabla_{\!\! xx} \bold{u}_{x,l})  
\end{bmatrix}}_{\bold{L}_{l}\in \mathbb{C}^{ n\times n}}
\underbrace{ \begin{bmatrix}
   \bold{m}_{\epsilon}  
\end{bmatrix}}_{\bold{m} \in \mathbb{R}^{n \times 1}}
=
 \underbrace{\begin{bmatrix}
    \omega^2 \bold{s}_l - \omega^2  \text{diag}(\bold{m}_{v_0})\bold{u}_{x,l} - \text{diag}(\bold{m}_\delta)\nabla_{\!\! zz} \bold{u}_{z,l}  
\end{bmatrix}
}_{\bold{y}_l \in \mathbb{C}^{n \times 1}}$ \Tstrut \\ [25pt] \hline
 p. & p. & a. & $\underbrace{
\begin{bmatrix}
    \text{diag}(\nabla_{\!\! zz} \bold{u}_{z,l}) \\
    \text{diag}(\nabla_{\!\! xx} \bold{u}_{x,l})    
\end{bmatrix}}_{\bold{L}_{l} \in \mathbb{C}^{2n \times n}}
\underbrace{ \begin{bmatrix}
   \bold{m}_{\delta} 
\end{bmatrix}}_{\bold{m} \in \mathbb{R}^ {n\times 1}}
=
 \underbrace{\begin{bmatrix}
    \omega^2 \bold{s}_l -\omega^2  \text{diag}(\bold{m}_{v_0})\bold{u}_{x,l} - \text{diag}(\bold{m}_\epsilon)\nabla_{\!\! xx} \bold{u}_{x,l}    \\
    \omega^2 \bold{s}_l -\nabla_{\!\! zz} \bold{u}_{z,l} - \omega^2  \text{diag}(\bold{m}_{v_0})\bold{u}_{z,l} 
\end{bmatrix}
}_{\bold{y}_l \in \mathbb{C}^{2n \times 1}}$ \Tstrut \\ [25pt] \hline
 p. & a. & a. & $\underbrace{
\begin{bmatrix}
     \text{diag}(\nabla_{\!\! xx} \bold{u}_{x,l})  & \text{diag}(\nabla_{\!\! zz} \bold{u}_{z,l}) \\
    \bold{0} &  \text{diag}(\nabla_{\!\! xx} \bold{u}_{x,l})    
\end{bmatrix}}_{\bold{L}_{l} \in \mathbb{C}^{2n \times 2n}}
\underbrace{ \begin{bmatrix}
   \bold{m}_{\epsilon}   \\
   \bold{m}_{\delta} 
\end{bmatrix}}_{\bold{m} \in \mathbb{R}^{2n \times 1}}
=
 \underbrace{\begin{bmatrix}
    \omega^2 \bold{s}_l-\omega^2  \text{diag}(\bold{m}_{v_0}) \bold{u}_{x,l}   \\
    \omega^2 \bold{s}_l -\nabla_{\!\! zz} \bold{u}_{z,l} - \omega^2  \text{diag}(\bold{m}_{v_0}) \bold{u}_{z,l}
\end{bmatrix}
}_{\bold{y}_l \in \mathbb{C}^{2n \times 1}}$ \Tstrut \\ [25pt] \hline
 a. & a. & p. & $\underbrace{
\begin{bmatrix}
    \omega^2  \text{diag}(\bold{u}_{x,l}) & \text{diag}(\nabla_{\!\! xx} \bold{u}_{x,l})   \\
   \omega^2  \text{diag}(\bold{u}_{z,l}) & \bold{0}    
\end{bmatrix}}_{\bold{L}_{l} \in \mathbb{C}^{2n \times 2n}}
\underbrace{ \begin{bmatrix}
   \bold{m}_{{v}_0}  \\
   \bold{m}_{\epsilon}    
\end{bmatrix}}_{\bold{m} \in \mathbb{R}^{2n \times 1}}
=
 \underbrace{\begin{bmatrix}
    \omega^2 \bold{s}_l - \text{diag}(\bold{m}_\delta)\nabla_{\!\! zz} \bold{u}_{z,l} \\
    \omega^2 \bold{s}_l -\nabla_{\!\! zz} \bold{u}_{z,l} -   \text{diag}(\bold{m}_\delta)\nabla_{\!\! xx} \bold{u}_{x,l}
\end{bmatrix}
}_{\bold{y}_l\in \mathbb{C}^{2n \times 1}}$ \Tstrut \\ [25pt] \hline
 a. & p. & a. & $\underbrace{
\begin{bmatrix}
   \omega^2 \text{diag}(\bold{u}_{x,l})   & \text{diag}(\nabla_{\!\! zz} \bold{u}_{z,l}) \\
   \omega^2  \text{diag}(\bold{u}_{z,l})  &  \text{diag}(\nabla_{\!\! xx} \bold{u}_{x,l})    
\end{bmatrix}}_{\bold{L}_{l} \in \mathbb{C}^{2n \times 2n}}
\underbrace{ \begin{bmatrix}
   \bold{m}_{v_0}  \\
   \bold{m}_{\delta} 
\end{bmatrix}}_{\bold{m} \in \mathbb{R}^{2n \times 1}}
=
 \underbrace{\begin{bmatrix}
    \omega^2 \bold{s}_l- \text{diag}(\bold{m}_\epsilon)\nabla_{\!\! xx} \bold{u}_{x,l}  \\
    \omega^2 \bold{s}_l -\nabla_{\!\! zz} \bold{u}_{z,l}
\end{bmatrix}
}_{\bold{y}_l \in \mathbb{C}^{2n \times 1}}$ \Tstrut \\ [25pt] \hline
\end{tabular}
\end{center}
\end{table*}
%
%
\subsection{ADMM-based acoustic VTI wavefield reconstruction inversion}
We consider the following bound-constrained TV-regularized multivariate optimization problem associated with the wave equation described by equation \ref{linearu}:
\begin{equation}   \label{main1_2}
\begin{aligned}
& \underset{\bold{U},\bold{m} \in \mathcal{C}}{\min}  \quad \sum  \sqrt{|\partial_x \bold{m}|^2 + |\partial_z\bold{m}|^2} ,\\
&\text{subject to} ~~~  \begin{cases} \bold{P}\bold{U} = \bold{D},\\ \bold{A(m)}\bold{U}=\bold{S}, \end{cases}
\end{aligned}
\end{equation} 
where $\partial_x$ and $\partial_z$ are, respectively, first-order finite-difference operators in the horizontal and vertical directions with appropriate boundary conditions, 
 $\mathcal{C} = \{\bold{x} \in \mathbb{R}^{3n\times 1}~\vert~ \bold{m}_{{l}} \leq \bold{x} \leq \bold{m}_{{u}}\}$ is the set of all feasible models bounded by a priori lower 
 bound $\bold{m}_{l}$ and upper bound $\bold{m}_{u}$,
 \begin{equation} 
\bold{D} =
\begin{bmatrix}
\bold{d}_{1} & \bold{d}_{2} & ... & \bold{d}_{n_s}
\end{bmatrix} \in \mathbb{C}^{n_r \times n_s} \nonumber
\end{equation}   
with $\bold{d}_l$ denoting the recorded data for the $l$th source, each including $n_r$ samples (the number of receivers), $\bold{P} \in \mathbb{R}^{n_r \times 2n}$ is a linear observation operator which samples the wavefields at the receiver positions.
Here, we assume that the sampling operator is identical across all sources (stationary-recording acquisitions). However, one may used a specific operator for each source.
 
We solve this constrained optimization problem with ADMM \citep{Boyd_2011_DOS,Aghamiry_2019_IBC}, an augmented Lagrangian method with operator splitting, leading to the following saddle point problem
\begin{equation}   \label{main1_3}
\begin{aligned}
& \underset{\bold{U},\bold{m} \in \mathcal{C}}{\min} ~ \underset{\bar{\bold{D}},\bar{\bold{S}} }{\max} \quad \sum  \sqrt{|\partial_x \bold{m}|^2 + |\partial_z\bold{m}|^2} \\
& + \bra \bar{\bold{D}} , \bold{P}\bold{U} - \bold{D} \ket + \lambda_0 \| \bold{P}\bold{U} - \bold{D} \|_{F}^2 \\
& + \bra \bar{\bold{S}}, \bold{A(m)}\bold{U} - \bold{S} \ket +  \lambda_1 \| \bold{A(m)}\bold{U} - \bold{S} \|_{F}^2,
\end{aligned}
\end{equation}
where $\| \bullet \|_{F}^2$ denotes the Frobenius norm of $\bullet$, $\lambda_0, \lambda_1  >0$ are penalty parameters and  $\bar{\bold{D}} \in \mathbb{C}^{n_r \times n_s} \nonumber$ and $\bar{\bold{S}} \in \mathbb{C}^{n \times n_s} \nonumber$ are the Lagrange multipliers (dual variables).

The last two lines of Equation \ref{main1_3} shows that the augmented Lagrangian function combines a Lagrangian function (left terms) and a penalty function (right terms).
Also, scaling the Lagrange multipliers by the penalty parameters, $\tilde{\bold{D}} = -\bar{\bold{D}}/\lambda_0$ and $\tilde{\bold{S}} = -\bar{\bold{S}}/\lambda_1$, allows us to recast the augmented Lagrangian function in a more convenient form \citep[][ Section 3.1.1]{Boyd_2011_DOS}
\begin{equation}   \label{main1_4}
\begin{aligned}
& \underset{\bold{U},\bold{m} \in \mathcal{C}}{\min} ~ \underset{\bar{\bold{D}},\bar{\bold{S}} }{\max} \quad \sum  \sqrt{|\partial_x \bold{m}|^2 + |\partial_z\bold{m}|^2} \\
& + \lambda_0 \| \bold{P} \bold{U} - \bold{D} - \tilde{\bold{D}} \|_{F}^2 - \| \tilde{\bold{D}} \|_{F}^2 \\
& +  \lambda_1 \| \bold{A(m)}\bold{U} - \bold{S} - \tilde{\bold{S}} \|_{F}^2 - \| \tilde{\bold{S}} \|_{F}^2,
\end{aligned}
\end{equation}
where the scaled dual variables have been injected in the penalty functions. 

In the WRI framework, the augmented Lagrangian method provides an efficient and easy-to-tune optimization scheme that extends the parameter search space by introducing a significant relaxation of the wave equation at the benefit of the observation equation during the early iterations, while satisfying the two equations at the convergence point.
We solve the saddle point problem, equation \ref{main1_4}, with the method of multiplier, in which the primal variables, $\bold{U}$ and $\bold{m}$, and the dual variables, $\tilde{\bold{D}}$ and $\tilde{\bold{S}}$, are updated in alternating mode. The dual problem is iteratively solved with basic gradient ascent steps. Accordingly, we immediately deduce from equation \ref{main1_4} that the scaled dual variables $\tilde{\bold{D}}$ and $\tilde{\bold{S}}$ are formed by the running sum of the constraint violations (the data and source residuals) in iterations. They update the RHSs (the data and the sources) of the quadratic penalty functions in equation \ref{main1_4} to refine the primal variables $\bold{U}$ and  $\bold{m}$ at a given iteration from the residual source and data errors (this RHS updating is a well known procedure to iteratively refine solutions of ill-posed linear inverse problems). The bi-variate primal problem, equation~\ref{main1_2}, is biconvex due to the bilinearity of the wave equation highlighted in the previous section. Therefore, it can be broken down into two linear subproblems for $\bold{U}$ and $\bold{m}$, which can be solved efficiently in alternating mode with ADMM after noting that the TV regularizer is convex \citep{Aghamiry_2019_IBC}. 

As pointed out by \citet{Aghamiry_2019_IWR}, a key advantage of augmented Lagrangian methods compared to penalty methods is that fixed penalty parameters $\lambda_0$ and $\lambda_1$ can be used in iterations, because the Lagrange multipliers progressively correct for the constraint violations generated by the penalty terms through the above mentioned RHS updating.   \\

Starting from an initial model $\bold{m}$ and zero-valued dual variables, the $k$th ADMM iteration embeds the following steps (see Appendix B for the complete development):
\subsubsection{Step 1: The primal wavefield reconstructions.}
Build regularized wavefields by solving the following multi-RHS system of linear equations with direct or iterative methods suitable for sparse matrices:
\begin{equation} \label{UE}
\begin{split} 
&\big[\lambda_0 \bold{P}^T\bold{P}+ \lambda_1 \bold{A}(\bold{m}^k)^T \bold{A}(\bold{m}^k)\big]\bold{U}=\\ &\big[\lambda_0 \bold{P}^T [\bold{D}+\tilde{\bold{D}}^k]
+\lambda_1\bold{A}(\bold{m}^k)^T [\bold{S}+\tilde{\bold{S}}^k]\big]. 
\end{split} 
\end{equation}
By choosing a small value of $\lambda_0/\lambda_1$, the reconstructed wavefields closely fit the observations during the early iterations, while only weakly satisfying the wave equation.
Problem \ref{UE} can be also interpreted as an extrapolation problem to reconstruct $\bold{U}$, when the observation equation (i.e. $\bold{PU=D}$) is augmented with the wave equation. 

To mitigate the computational burden of the wavefield reconstruction, we solve equation \ref{UE} with a fourth-order wave equation operator following the parsimonious approach of \citet{Operto_2014_FAT}, while the subsequent model estimation subproblem 
relies on the bilinear wave equation provided in equation \ref{sovtit}. The elimination procedure allowing to transform the system of two second-order wave equations for $\bold{u}_x$ and $\bold{u}_z$, equation \ref{sovtit}, into a fourth-order wave
equation for $\bold{u}_x$ coupled with the closed-form expression of $\bold{u}_z$ as a function of $\bold{u}_x$ is reviewed in Appendix C.

\subsubsection{Step 2: The primal model estimation.}
\vspace{-0.1cm}
The reconstructed wavefields, equation \ref{UE}, are injected in the linearized equation \ref{linearma} to update the subsurface parameters by solving the following system of linear equations:
\begin{equation} \label{ME}
\begin{split}
&\left[\lambda_1\sum_{l=1}^{n_s}\bold{L}_l^T\bold{L}_l + \overline{\nabla}^T\Gamma\overline{\nabla} + {Z}\bold{I}\right]\bold{m}= 
\lambda_1\sum_{l=1}^{n_s}\bold{L}_l^T(\bold{y}^k_l+\tilde{\bold{s}}^k_l)\\
&\hspace{3cm} + \overline{\nabla}^T\Gamma(\bold{p}^k+\tilde{\bold{p}}^k) 
+ {Z}(\bold{q}^k+\tilde{\bold{q}}^k),
\end{split}
\end{equation}
where 
\begin{equation} \label{partial}
\overline{\nabla}=
\begin{bmatrix}
\partial_x \\
\partial_z
\end{bmatrix}.
\end{equation}
In equation \ref{ME}, $\bold{p}$ and $\bold{q}$ are auxiliary primal variables, which have been introduced to solve the bound-constrained TV-regularized parameter estimation subproblem with the split Bregman method 
(Appendix B).
The vectors $\tilde{\bold{p}}$ and $\tilde{\bold{q}}$ are the corresponding dual variables. These auxiliary primal and dual variables are initialized to $\bold{0}$ during the first ADMM iteration.

The operator $Z$ is a diagonal weighting matrix defined as
 \begin{equation} \label{Zeta}
 Z=
 \begin{bmatrix}
 \zeta_{v_0} \bold{I} & \bold{0} & \bold{0}\\
 \bold{0} & \zeta_{\epsilon}\bold{I} & \bold{0}\\
  \bold{0} & \bold{0} & \zeta_{\delta}\bold{I}\\
\end{bmatrix} \in \mathbb{R}_+^{3n \times 3n}, 
 \end{equation}
where $\zeta_{v_0}, \zeta_{\epsilon}, \zeta_{\delta}>0$ control the relative weights assigned to the bound constraints applied on the three parameter classes. Note that the bound constraints introduce also a damping (DMP) or zero-order
Tikhonov regularization in the Hessian of equation \ref{ME}.

In the same way,
  $\bold{\Gamma}$ is a diagonal matrix defined as
 \begin{equation} \label{Gamma}
 \bold{\Gamma}=\begin{bmatrix}
 \bold{\Gamma}_{11} & \bold{0} \\
 \bold{0} & \bold{\Gamma}_{22} \\
\end{bmatrix},
 \end{equation}
where
 \begin{equation}
 \bold{\Gamma}_{11}= \bold{\Gamma}_{22}=
 \begin{bmatrix}
 \gamma_{v_0}\bold{I} & \bold{0} & \bold{0}\\
 \bold{0} & \gamma_{\epsilon}\bold{I} & \bold{0}\\
  \bold{0} & \bold{0} & \gamma_{\delta}\bold{I}\\
\end{bmatrix} \in \mathbb{R}_+^{3n \times 3n}, \nonumber
 \end{equation}
and $\gamma_{v_0}, \gamma_{\epsilon}, \gamma_{\delta}>0$ control the soft thresholding that is performed by the TV regularizer, equation \ref{prox0}.  We remind that augmented Lagrangian methods
seek to strictly satisfy the constraints at the convergence point only, not at each iteration. Therefore, the relative values of these penalty parameters have a significant impact upon the path followed by the inversion to converge toward this convergence point.\\
\subsubsection{Step 3: The TV primal update.}
Update the TV primal variable $\bold{p}$ via a TV proximity operator.
Set
\begin{equation}
\bold{z} \leftarrow \overline{\nabla}\bold{m}^{k+1}-\tilde{\bold{p}}^k=
\begin{bmatrix}
\bold{z}_x\\
\bold{z}_z
\end{bmatrix}, \nonumber
\end{equation}
then 
\begin{equation} \label{primal_p}
\bold{p}^{k+1}\leftarrow\text{prox}_{\Gamma^{-1}}(\bold{z})=
\begin{bmatrix}
\xi\circ \bold{z}_x\\
\xi\circ \bold{z}_z
\end{bmatrix},
\end{equation}
where
\begin{equation} \label{prox0}
\xi = \max(1 - \frac{1}{\bold{\Gamma} \sqrt{\bold{z}_x^2+\bold{z}_z^2}},0)
\end{equation}  
and $\bold{x}\circ \bold{y}$ denotes the Hadamard (component-wise) product of $\bold{x}$ and $\bold{y}$. Also, the power of 2 indicates the Hadamard product of $\bold{x}$ with itself, i.e. $\bold{x}^2 =\bold{x} \circ \bold{x}$. 
\subsubsection{Step 4: The bounding constraint primal update.}
Update the primal variable $\bold{q}$ via a projection operator, which has the following component-wise form
 \begin{equation} \label{dual_q}
 \bold{q}^{k+1}\leftarrow
\text{proj}_{\mathcal{C}} (\bold{m}^{k+1} - \tilde{\bold{q}}^k),
 \end{equation}
where the projection operator is given by
\begin{equation}
\text{proj}_{\mathcal{C}} (\bold{x}) = \min(\max(\bold{x},\bold{m}_{l}), \bold{m}_{u}). \nonumber
\end{equation}

\subsubsection{Step 5: Dual updates.}
Update the scaled dual variables with gradient ascent steps
\begin{equation}
\begin{cases}
\tilde{\bold{S}}^{k+1} &\leftarrow  \tilde{\bold{S}}^{k}  +\bold{S}- \bold{A}(\bold{m}^{k+1})\bold{U}^{k+1} , \\
\tilde{\bold{D}}^{k+1} &\leftarrow  \tilde{\bold{D}}^{k} + \bold{D}- \bold{P} \bold{U}^{k+1}, \\
\tilde{\bold{p}}^{k+1} &\leftarrow \tilde{\bold{p}}^{k} + \bold{p}^{k+1}-\overline{\nabla} \bold{m}^{k+1},\\
\tilde{\bold{q}}^{k+1} &\leftarrow \tilde{\bold{q}}^{k} + \bold{q}^{k+1}- \bold{m}^{k+1},
\end{cases}
\end{equation}

\subsubsection{Step 6: Check the stopping condition.}
Exit if the preset stopping conditions are satisfied else go to step 1. We will describe the stopping criteria of iterations in the following "Numerical example" section for each numerical example

\subsection{Hyperparameter tuning}
We tune the different penalty parameters by extending the procedure reviewed by \citet{Aghamiry_2019_IBC,Aghamiry_2019_CROb} to multiparameter reconstruction.

We start from the last subproblem of the splitting procedure and set the penalty parameters contained in $\bold{\Gamma}$. These hyperparameters control the soft thresholding performed by the TV regularization, equation~\ref{primal_p}. We set $\gamma_i=2\% \max {\sqrt{\bold{z}^2_{i_x} + \bold{z}^2_{i_z}}}$, where the subscript $i \in \{v_0,\epsilon,\delta \}$ denotes the parameter class ($\bold{m}_{v_0}$, $\bold{m}_\epsilon$, $\bold{m}_\delta$). This tuning can be refined by using a different thresholding percentage for each parameter class adaptively during iterations or according to prior knowledge of the geological structure, coming  from well logs for example. Also, we use the same weight for the damping regularization associated with the bound constraints and the TV regularization: $\zeta_i = \gamma_i$.  \\
Then, we select $\lambda_1$ as a percentage of the mean absolute value of the diagonal coefficients of $\sum_{l=1}^{n_s}\bold{L}_l^T\bold{L}_l$ during the parameter estimation subproblem, equation \ref{ME}. This percentage is set according to the weight that we want to assign to the TV regularization and the bound constraints relative to the wave equation constraint during the parameter estimation. Parameter $\lambda_1$ may be increased during iterations to reduce the weight of TV regularization and bound constraints near the convergence point. We found this adaptation useful when we start from very crude initial models. Finally, we set $\lambda_0$ such that  $\lambda=\lambda_1/\lambda_0$ is a small fraction of the highest eigenvalue $\xi$ of the normal operator $\bold{A(m})^{-T}\bold{P}^{T}\bold{P}\bold{A(m})^{-1}$ during the wavefield reconstruction subproblem, equation~\ref{UE}, according to the criterion proposed by \citet{vanLeeuwen_2016_PMP}. In all the numerical tests, we use  $\lambda=1\text{e-2} \xi$ and $\lambda=1\text{e-0} \xi$ for noiseless and noisy data, respectively. This tuning of $\lambda$ is indeed important because it controls the extension of the search space. 
A too high value of $\lambda$ reduces the weight of $\|\bold{PU-D}\|_2^2$ during the wavefield reconstruction and makes IR-WRI behave like a reduced approach. Conversely, using a small value for $\lambda$ fosters data fitting and expends the search space accordingly. However, a too small value can lead to a prohibitively high number of iterations of the augmented Lagrangian method before the wave equation constraint is fulfilled with sufficient accuracy. Moreover, when data are contaminated by noise, a too small value for $\lambda$ will make the wavefield reconstruction over-fit the data and drive the algorithm to be a poor minimizer.
We always use $\lambda$ as a fixed percentage of $\xi$ in iterations.\\
\section{Numerical examples}
%
%
\subsection{Inclusion test}
We first assess multi-parameter IR-WRI with a simple inclusion example for noiseless data. The experimental setup in terms of model, acquisition geometry and frequency selection is identical to that used by \citet{Gholami_2013_WPA1}. The vertical velocity $\bold{v}_0$ and the Thomsen's parameters $\bold{\delta}$ and $\bold{\epsilon}$ are 3 km/s, 0.05 and 0.05, respectively, in the homogeneous background model and 3.3 km/s, 0.1 and 0.2, respectively, in the inclusion of radius 100 m. Nine frequencies between 4.8 Hz and 19.5 Hz are processed simultaneously during IR-WRI and a maximum of 25 iterations is used as stopping criterion for iterations. An ideal fixed-spread acquisition is used, where 64 sources and 320 receivers surround the inclusion, hence providing a complete angular illumination of the target. 

Although we use $(1/\bold{v}_0^2,\sqrt{1 + 2 \delta},1 + 2\epsilon)$ as optimization parameters during our inversions, we show the reconstructed models under the form of $\bold{v}_0$, $\delta$ and $\epsilon$ for comparison with the results of \citet{Gholami_2013_WPA1}. Note that the radiation patterns of the $(1/\bold{v}_0^2,\sqrt{1 + 2 \delta},1 + 2\epsilon)$ parametrisation are scaled versions of those of the $(\bold{v}_0,\delta,\epsilon)$ parametrisation: namely, they exhibit the same amplitude variation with scattering angle, with however different amplitudes. This means that, for an equivalent regularization and parameter scaling, we expect similar resolution and trade-off effects with these two parametrisations. Note that $\bold{v}_0$ is provided in km/s in our inversion such that the order of magnitude of $1/\bold{v}_0^2$ is of the order of $\delta$ and $\epsilon$.  

\begin{figure}
\includegraphics[scale=0.45]{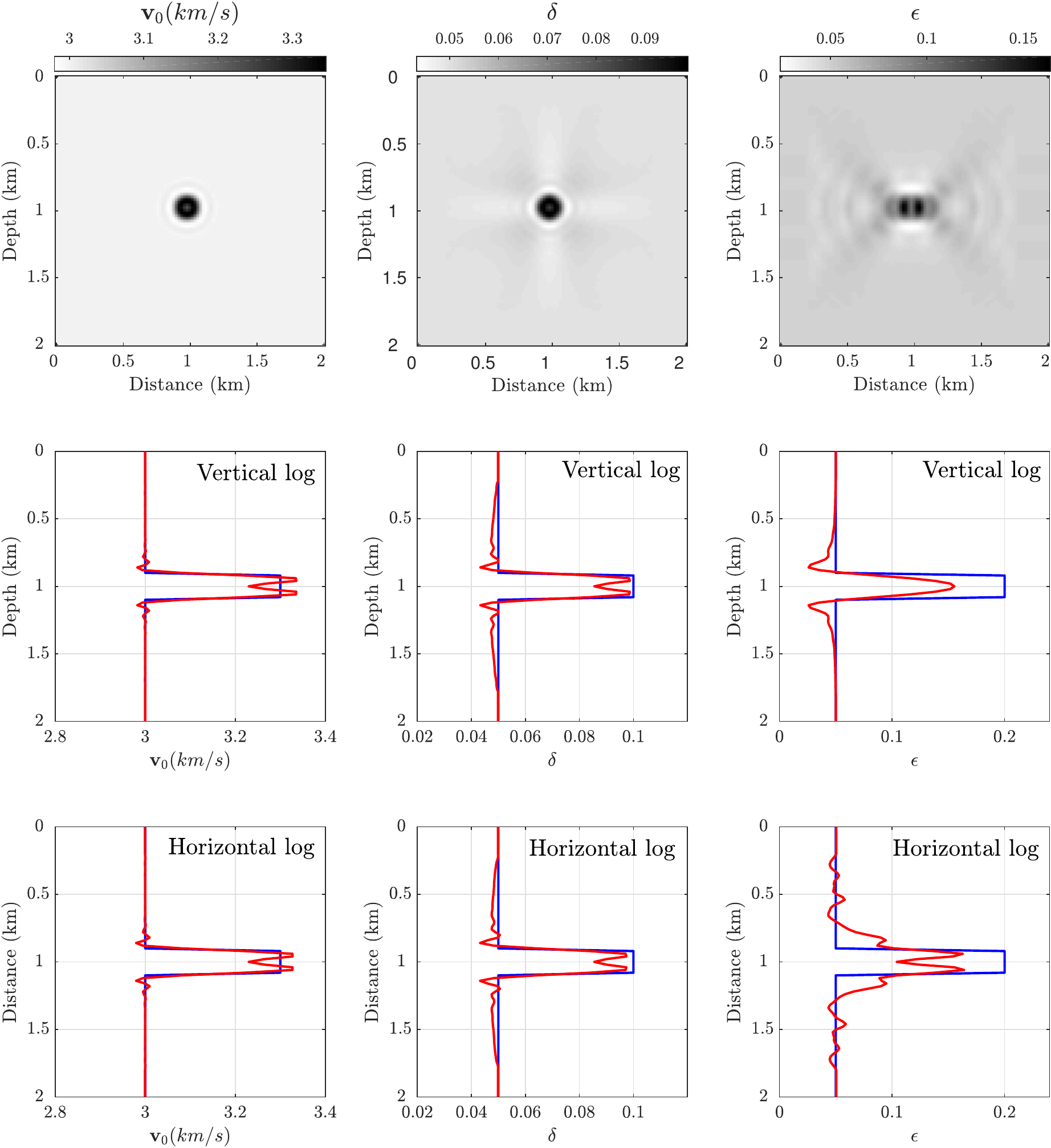}
\caption{Inclusion test: Mono-parameter IR-WRI results.  For $\bold{v}_0$ (first column), $\bold{\delta}$ (second column), and $\bold{\epsilon}$ (third column). The initial models are the true homogeneous background models. The vertical and horizontal profiles in the true model (blue) and estimated model (red) are extracted across the center of the inclusions.}
\label{mono}
\end{figure}

\begin{figure}
\includegraphics[scale=0.45]{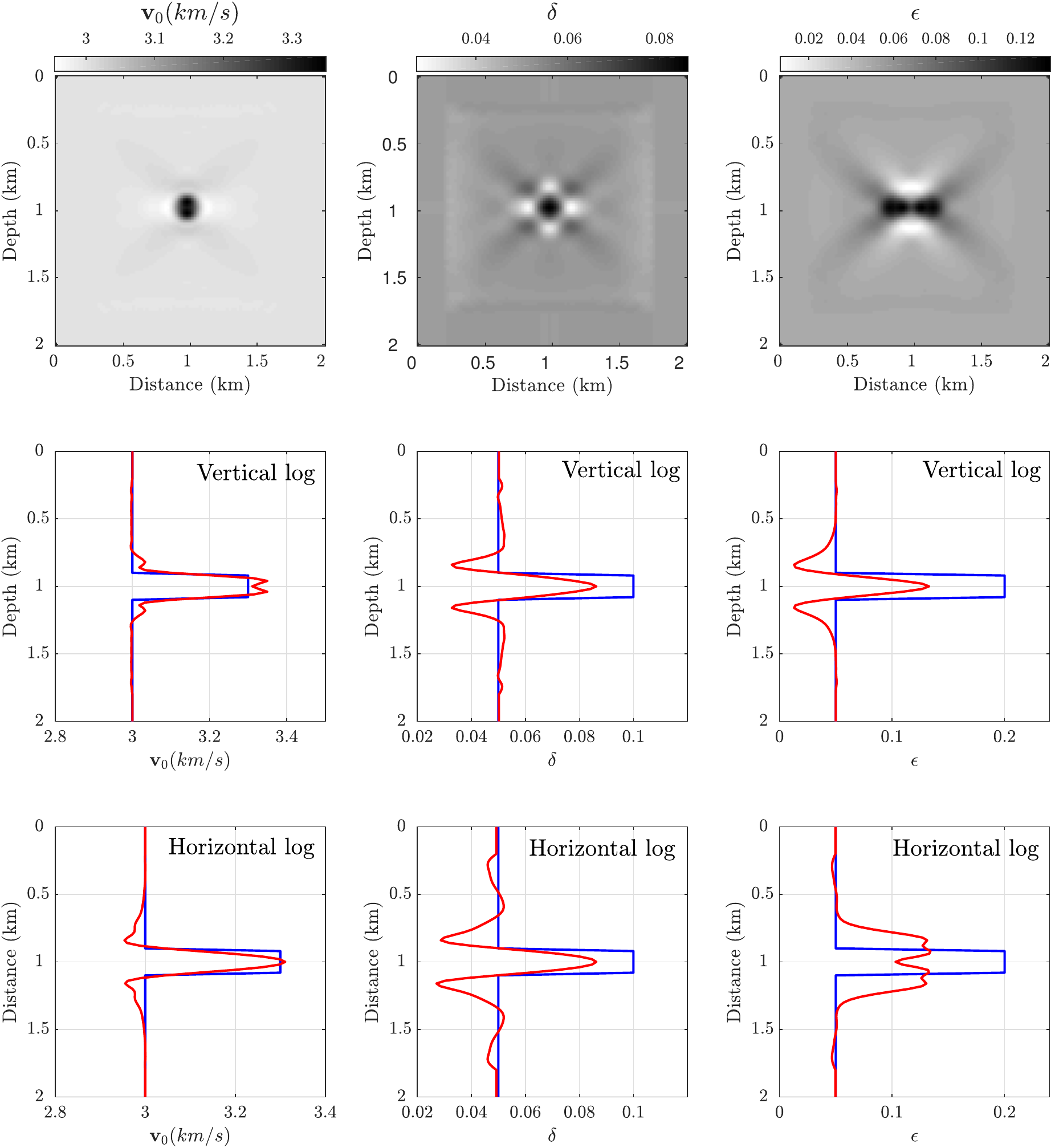}
\caption{Inclusion test: Multiparameter IR-WRI results for joint update of $\bold{v}$ (first column), $\bold{\delta}$ (second column), and $\bold{\epsilon}$ (third column). The initial models are the homogeneous background models.}
\label{joint}
\end{figure}

We start with bound-constrained IR-WRI with damping (DMP) regularization only ($\gamma_i=0, i \in \{v_0,\epsilon,\delta\}$ in equation \ref{ME}) and perform three independent mono-parameter reconstructions for $\bold{m}_{v_0}$, $\bold{m}_\delta$ and $\bold{m}_\epsilon$, respectively. For each mono-parameter inversion, the true model associated with the optimization parameter contains the inclusion, while the true models associated with the two passive parameters are homogeneous (Fig. \ref{mono}). 
For all three inversions, the starting models are the true homogeneous background models.
In other words, the data residuals contain only the footprint of the mono-parameter inclusion to be reconstructed. This test is used to assess the intrinsic resolution of IR-WRI for each parameter reconstruction, independently from the cross-talk issue \citep{Gholami_2013_WPA1}. It is reminded that this intrinsic resolution is controlled by the frequency bandwidth, the angular illumination provided by the acquisition geometry and the radiation pattern of the optimization parameter in the chosen subsurface parametrisation. As this test fits the linear regime of classical FWI, we obtain results (Fig. \ref{mono}) very similar to those obtained by \citet[][ Their Fig. 9]{Gholami_2013_WPA1}, where the shape of the reconstructed inclusions is controlled by the radiation pattern of the associated parameter. The reader is referred to \citet{Gholami_2013_WPA1} for a detailed analysis of these radiation patterns.  \\
Then, we perform the joint reconstruction of $\bold{m}_{v_0}$, $\bold{m}_{\delta}$ and $\bold{m}_{\epsilon}$, when the true model contains an inclusion for each parameter class (Fig. \ref{joint}). With our parameter scaling and subsurface parametrisation, $\bold{v}_0$, $\bold{\epsilon}$ and $\bold{\delta}$ are reconstructed with well balanced amplitudes compared to the results of \citet[][ Their Fig. 10]{Gholami_2013_WPA1}, where $\bold{v}_0$ has a dominant imprint in the inversion. Comparing the models reconstructed by the mono-parameter and multi-parameter inversions highlights however the wavenumber leakage generated by parameter cross-talks (Figs. \ref{mono} and \ref{joint}).

\begin{figure}
\includegraphics[scale=0.45]{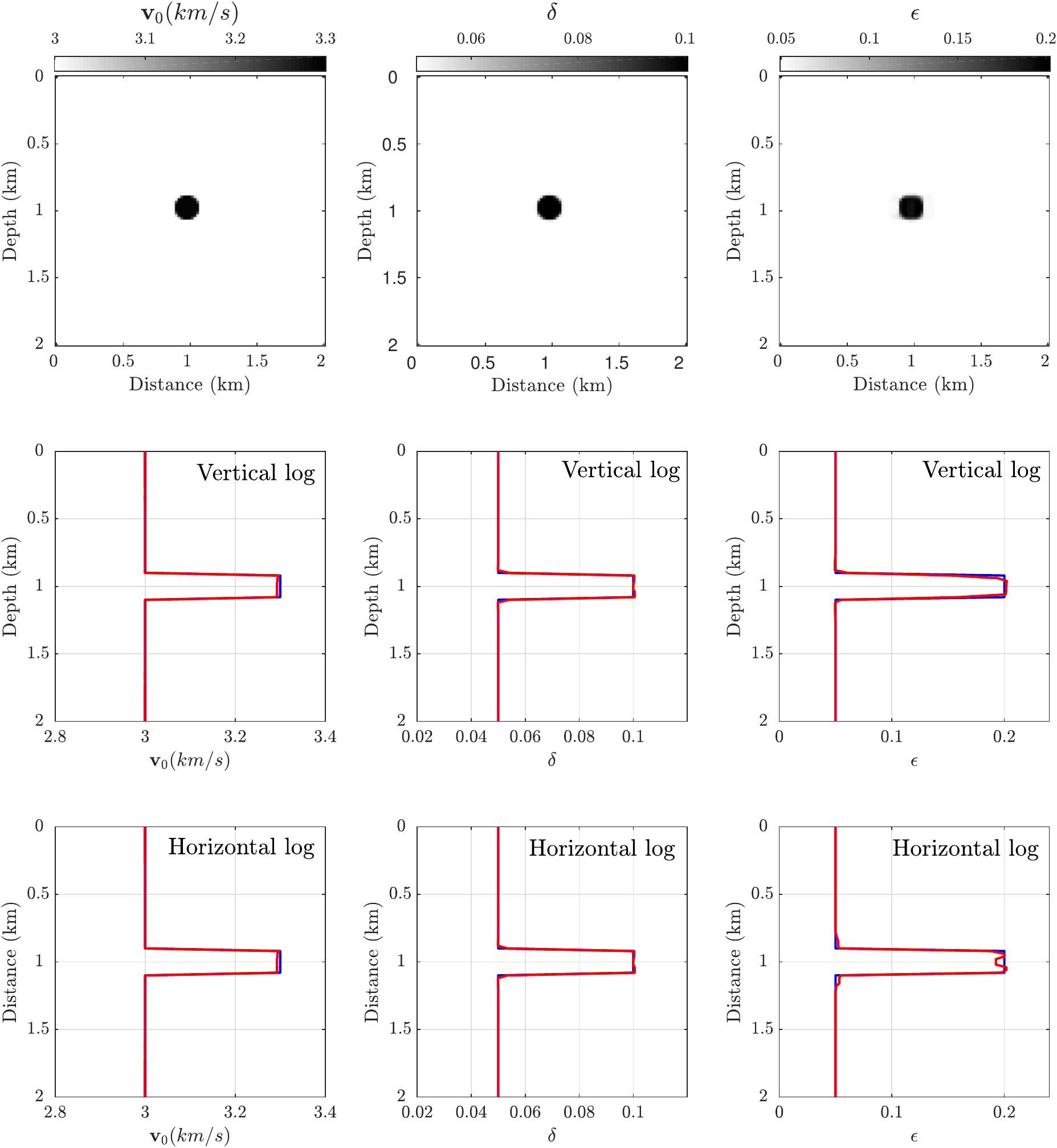}
\caption{Inclusion test: Same as Fig. \ref{mono}, but with TV regularization.}
\label{mono_TV}
\end{figure}

\begin{figure}
\includegraphics[scale=0.45]{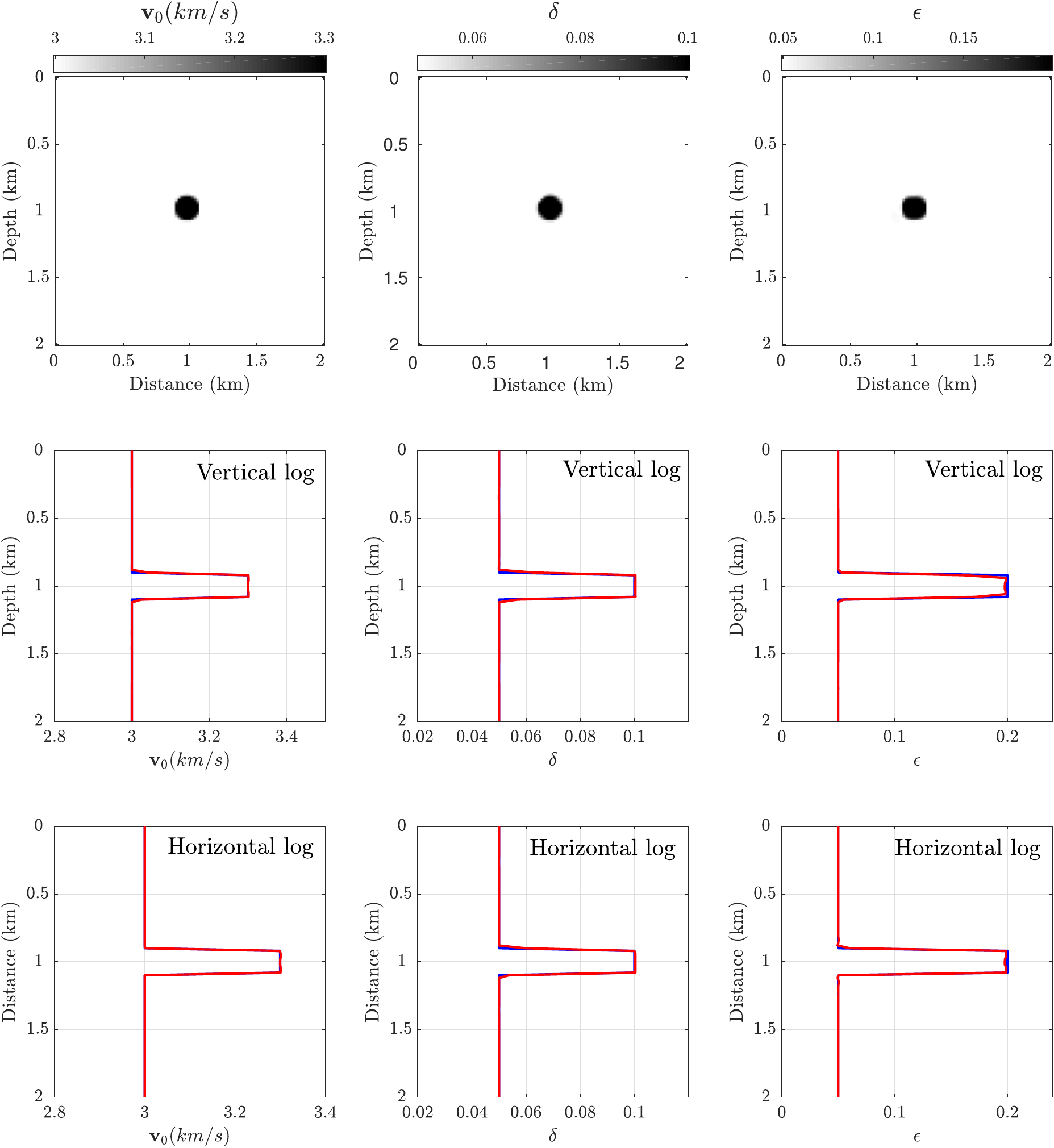}
\caption{Inclusion test: Same as Fig. \ref{joint}, but with TV regularization.}
\label{joint_TV}
\end{figure}

%
Then, we complement DMP regularization with TV regularization during bound-constrained IR-WRI for the above mono-parameter and multi-parameter experiments (Figs. \ref{mono_TV} and \ref{joint_TV}). The results show how TV regularization contributes to remove the wavenumber filtering performed by radiation patterns (compare for example Figs. \ref{mono} and \ref{mono_TV}) as well as the cross-talk artifacts during the multi-parameter inversion (compare Figs. \ref{joint} and \ref{joint_TV}).
Although this toy example has been designed with a piecewise constant model which is well suited for TV regularization, yet it highlights the potential role of TV regularization to mitigate the ill-posedness of FWI resulting from incomplete wavenumber coverage and parameter trade-off.
%
\subsection{Synthetic North Sea case study}
\subsubsection{Experimental setup}
\begin{figure*}
\centering
\includegraphics[width=0.8\textwidth]{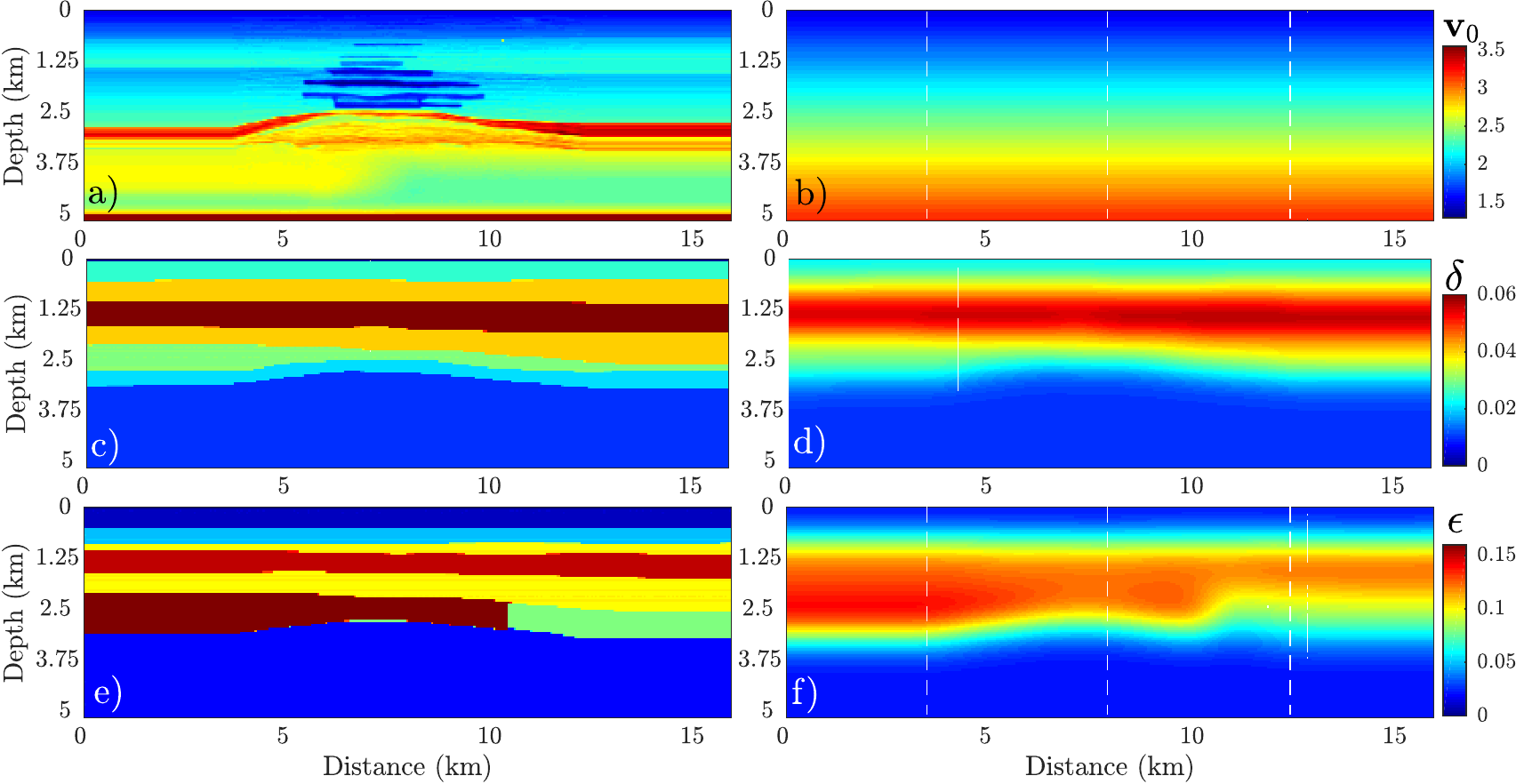} 
\caption{North Sea case study. (a) True $\bold{v}_0$ model. (b) Initial $\bold{v}_0$ model. (c) True $\delta$ model. (d)  Smoothed $\delta$ model, which is used as passive background model during inversion tests. (e) True $\epsilon$ model. (f) Smoothed $\epsilon$ model, which is used as a passive background model during the mono-parameter inversion and as an initial model during the joint reconstruction of $\bold{v}_0$ and $\bold{\epsilon}$. The vertical dashed lines in (b) and (f) indicate the location of vertical logs.}
\label{fig:val_initial}
\end{figure*}

We consider now a more realistic $16~km \times 5.2~km$ shallow-water model representative of the North Sea \citep{Munns_1985_VFG}. The reader is also referred to \citet{Gholami_2013_WPA2} for an application of acoustic VTI FWI on this model. The true model and the initial (starting) models for $\bold{v}_0$, $\bold{\delta}$ and $\bold{\epsilon}$ are shown in Fig. \ref{fig:val_initial}.
The subsurface model is formed by soft sediments in the upper part, a pile of low-velocity gas layers above a chalk reservoir, the top of which is indicated by a sharp positive velocity contrast at around 2.5~km depth, and a flat reflector at 5~km depth (Fig. \ref{fig:val_initial}a).
The initial $\bold{v}_0$ model is laterally homogeneous with velocity linearly increasing with depth between 1.5 to 3.2~km/s (Fig. \ref{fig:val_initial}b), while the  $\bold{\delta}$ and $\bold{\epsilon}$ initial models are Gaussian filtered version of the true models. Note that our initial/background  $\bold{v}_0$, $\bold{\delta}$ and $\epsilon$ models are cruder than those used by \citet{Gholami_2013_WPA2}.
%
%
%
%
The fixed-spread surface acquisition consists of 80 (reciprocal) explosive sources spaced 200~m apart at 75~m depth on the sea bottom and 320 (reciprocal) hydrophone receivers spaced 50~m apart at 25~m depth. Accordingly, the pressure wavefield is considered for the acoustic VTI inversion.
A free-surface boundary condition is used on top of the grid and the source signature is a Ricker wavelet with a 10~Hz dominant frequency.
We compute the recorded data in the true models (Figs. \ref{fig:val_initial}a,c,e) with the forward modelling engine described in Appendix C. During IR-WRI, we use the same forward engine to compute the modelled data according to an inverse crime procedure.
%
%
Common-shot gathers computed in the true model and in the initial model for a shot located at 14~km are shown in Fig. \ref{fig:val_seismo_true}. The seismograms computed in the true model are dominated by the direct wave, the diving waves from the sedimentary overburden, complex packages of pre- and post-critical reflections from the gas layers, the top of the reservoir and the deep reflector. The refracted wave from the deep reflector is recorded at a secondary arrival between -14km and -10km offset in Fig.  \ref{fig:val_seismo_true}a. Also, energetic reverberating P-wave reflections are generated by the wave guide formed by the shallow water layer and the weathering layer. They take the form of leaking modes with phase velocities higher than the water wave speed \citep{Operto_2018_MFF}. The seismograms computed in the starting model mainly show the direct wave and the diving waves, these latter being highly cycle skipped relative to those computed in the true model. 

%
%
We perform both mono-parameter IR-WRI for $\bold{m}_{v_0}$ and  multi-parameter IR-WRI for $\bold{m}_{v_0}$ and $\bold{m}_{\epsilon}$. In both cases, we compare the results that are obtained when bound-constrained IR-WRI is performed with DMP regularization only ($\gamma_i=0$) and with DMP + TV regularization. When $\bold{m}_{\epsilon}$ is involved as an optimization parameter, we had to introduce an additional regularization term $\| \bold{m} - \bold{m}_\epsilon^0\|$ in the parameter-estimation subproblem, equation \ref{main4_m}, in order to force the updates of $\epsilon$ to be smooth and close to $\epsilon^0$. The background models used for the passive anisotropic parameters (either $\epsilon$ and $\delta$ or $\delta$ alone) are the smooth models (Fig. \ref{fig:val_initial}d,f), which means that the IR-WRI results will be impacted upon by the smoothness of the passive parameters.

We perform IR-WRI with small batches of two frequencies with one frequency overlap between two consecutive batches, moving from the low frequencies to the higher ones according to a classical frequency continuation strategy. The starting and final frequencies are 3~Hz and 15~Hz and the sampling interval in one batch is 0.5~Hz. The stopping criterion for iterations and for each batch is given by
$k_{max}=15$ or
\begin{eqnarray}
\label{Stop}
\| \bold{A(m}^{k+1})\bold{U}^{k+1}-\bold{S}\|_{F} \leq \varepsilon_b~\& ~  \|\bold{PU}^{k+1}-\bold{D}\|_{F} \leq \varepsilon_d,
\end{eqnarray} 
where $k_{max}$ denotes the maximum iteration count, $\varepsilon_b$=1e-3, and $\varepsilon_d$=1e-5 for noiseless data and $\varepsilon_b$=1e-3 , $\varepsilon_d$= noise level of batch for noisy data.  

We perform three paths through the frequency batches to improve the IR-WRI results, using the final model of one path as the initial model of the next one (these cycles can be viewed as outer iterations of IR-WRI). The starting and finishing frequencies of the paths are [3, 6], [4, 8.5], [6, 15]~Hz respectively, where the first element of each pair shows the starting frequency and the second one is the finishing frequency.
\subsubsection{Convexity and sensitivity analysis}
\begin{figure}[h]
\centering
\includegraphics[width=0.5\textwidth]{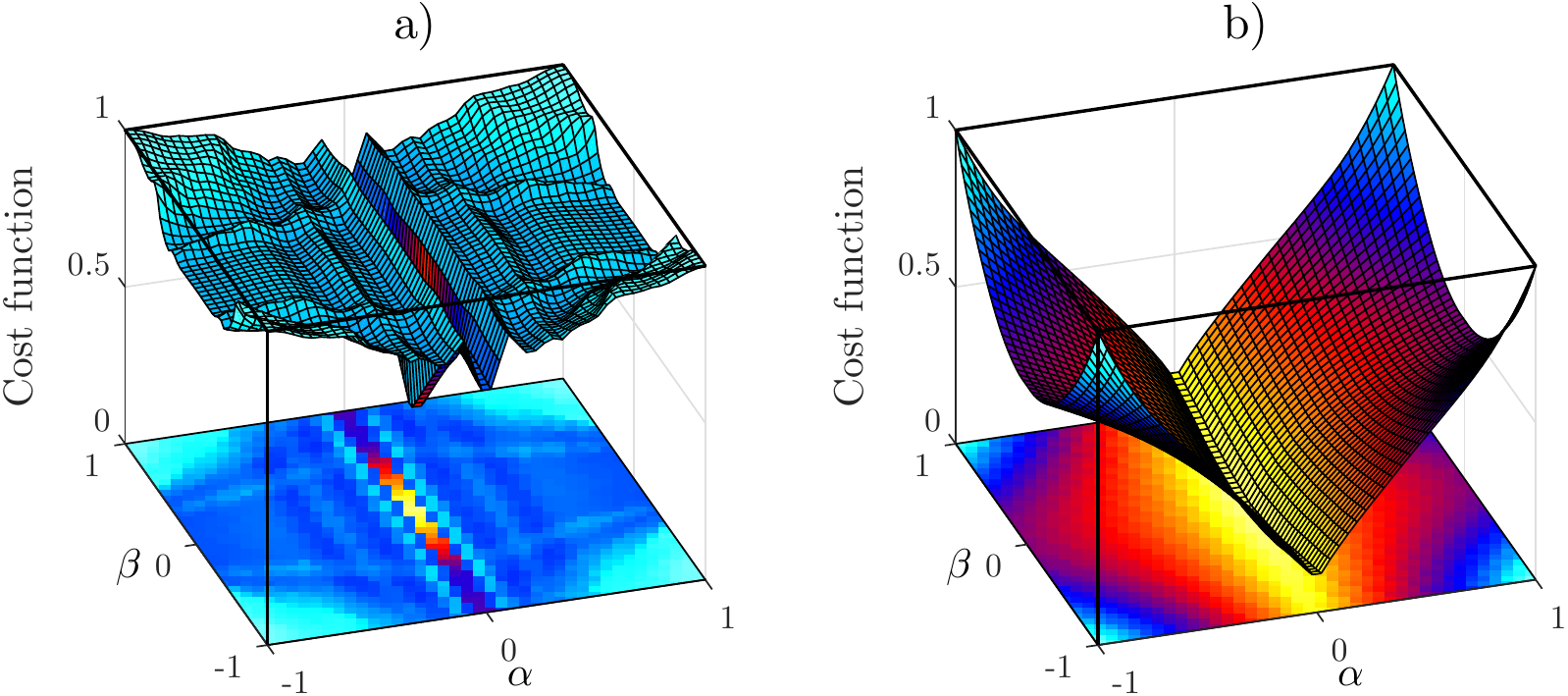}
\caption{North Sea case study. 3-Hz misfit function for the first iteration when $\bold{v}_0$ and $\bold{\epsilon}$ are active parameters and $\bold{\delta}$ is kept fixed to the smooth background model (Fig. \ref{fig:val_initial}f). (a) Classical reduced approach, (b) WRI. The variables $\alpha$ and $\beta$ parametrize the $\bold{v}_0$ and $\bold{\epsilon}$ models, respectively, for which the misfit function is computed (see text for details).}
\label{fig:val_cost}
\end{figure}

\begin{figure}
\centering
   \includegraphics[width=0.5\textwidth,clip=true,trim=0cm 0cm 0cm 0cm]{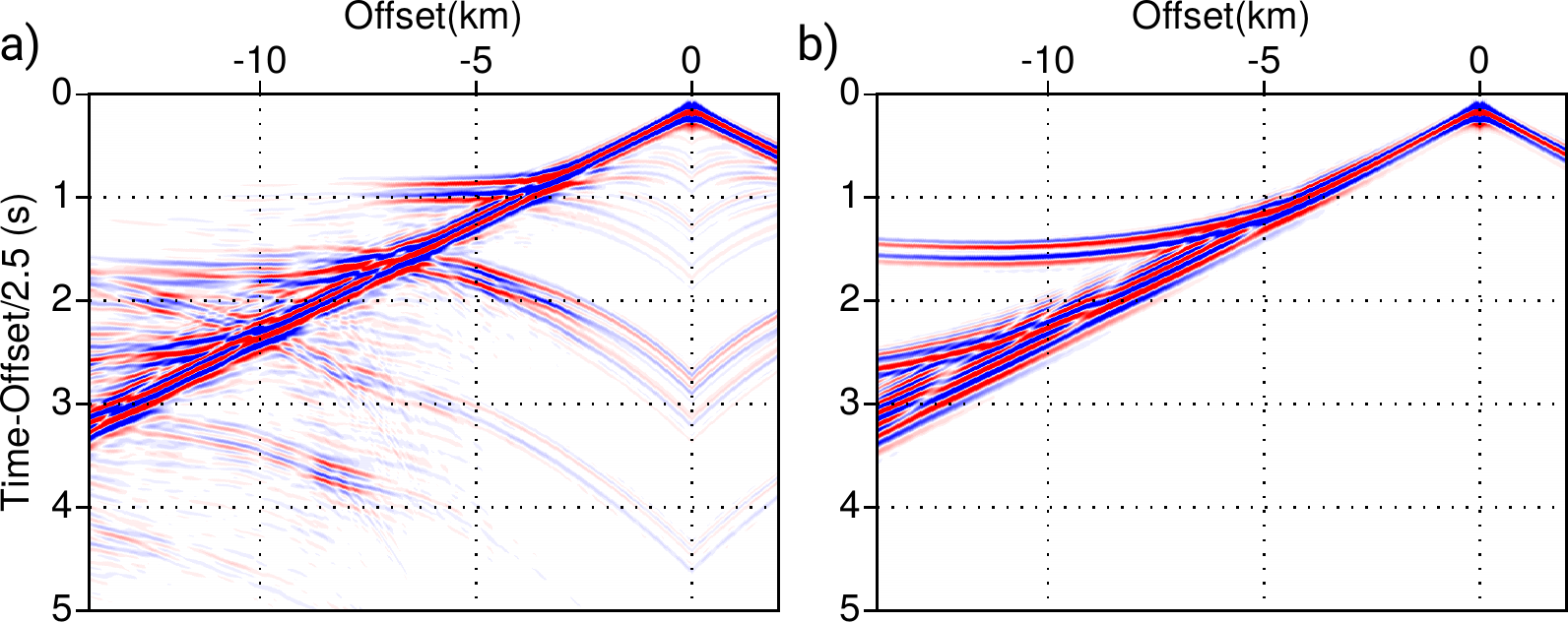} 
\caption{North Sea case study. Time domain seismograms computed in (a) the true model (Fig. \ref{fig:val_initial}a,c,e) and in (b) the initial model (Fig. \ref{fig:val_initial}b,d,f). The seismograms are plotted with a reduction velocity of 2.5 km/s.}
\label{fig:val_seismo_true}
\end{figure}

Before discussing the IR-WRI results, we illustrate how WRI extends the search space of FWI for the North Sea case study. For this purpose,  we compare the shape of the FWI misfit function with that of the parameter-estimation WRI subproblem for the 3~Hz frequency and for a series of $\bold{v}_0$ and $\epsilon$ models that are generated according to 
\begin{subequations}
\begin{eqnarray}
\bold{v}_0(\alpha)&=&\bold{v}_{0_{true}} + |\alpha| [\bold{v}_{{0}_{init}}-\bold{v}_{0_{true}}], \\
\bold{\epsilon}(\beta) &=&\bold{\epsilon}_{true} + |\beta| [\bold{\epsilon}_{init}-\bold{\epsilon}_{true}],
\end{eqnarray}
\end{subequations}
where $\bold{v}_{0_{true}}$ and $\bold{v}_{{0}_{init}}$ denote the true and the initial $\bold{v}_{0}$ models, respectively (Fig. \ref{fig:val_initial}a,b) and $-1 \leq \alpha \leq 1$. Similarly, $\bold{\epsilon}_{true}$ and $\bold{\epsilon}_{init}$ are the true and the initial $\bold{\epsilon}$ models, respectively (Fig. \ref{fig:val_initial}e,f) and $-1 \leq \beta \leq 1$. Finally, we use the true $\bold{\delta}$ model (Fig. \ref{fig:val_initial}c) to generate the recorded data and the smoothed version (Fig. \ref{fig:val_initial}d) as a passive parameter to evaluate the misfit function. The misfit functions of the classical reduced-approach FWI as well as that of WRI are shown in Fig. \ref{fig:val_cost}. 
The WRI misfit function is convex, while that of FWI exhibits spurious local minima along both the $\alpha$ and $\beta$ dimensions. Also, the sensitivity of the misfit function to $\bold{v}_0$ is much higher than that of $\bold{\epsilon}$ for the considered range of models as already pointed out by \citet{Gholami_2013_WPA1}, \citet{Gholami_2013_WPA2} and \citet[][ Their Fig. 2]{Cheng_2016_MEA}. The weaker sensitivity of the misfit function to $\epsilon$ makes the joint update of  $\bold{v}_{0}$ and  $\bold{\epsilon}$ challenging. \\
%
%
\subsubsection{Mono-parameter IR-WRI}
\begin{figure}[h]
\centering
   \includegraphics[width=0.46\textwidth]{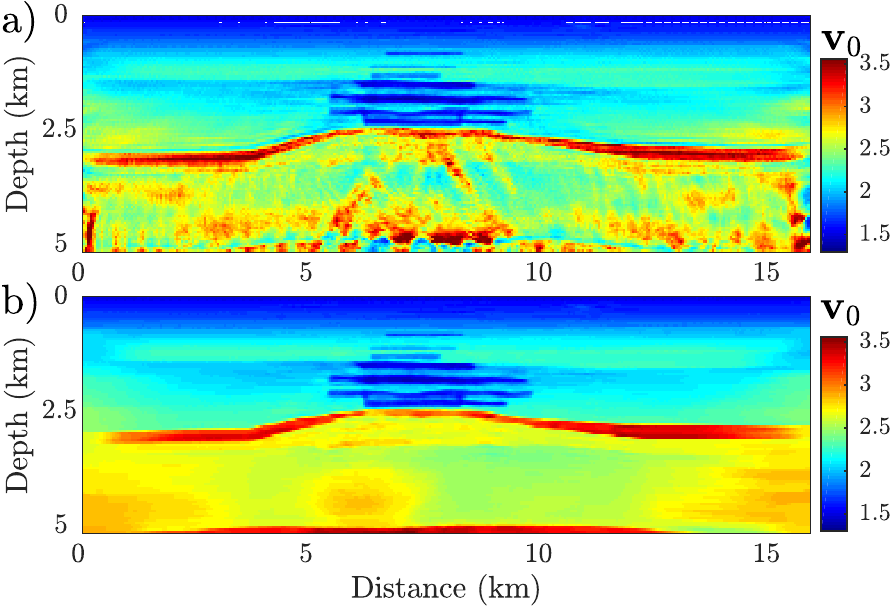} 
\caption{North Sea case study with noiseless data. Mono-parameter IR-WRI. Vertical wavespeed $\bold{v}_0$ models inferred from IR-WRI with (a) DMP regularization, (b) TV regularization.}
\label{fig:val_mono_noiseless}
\end{figure}

\begin{figure}
\centering
   \includegraphics[width=0.48\textwidth]{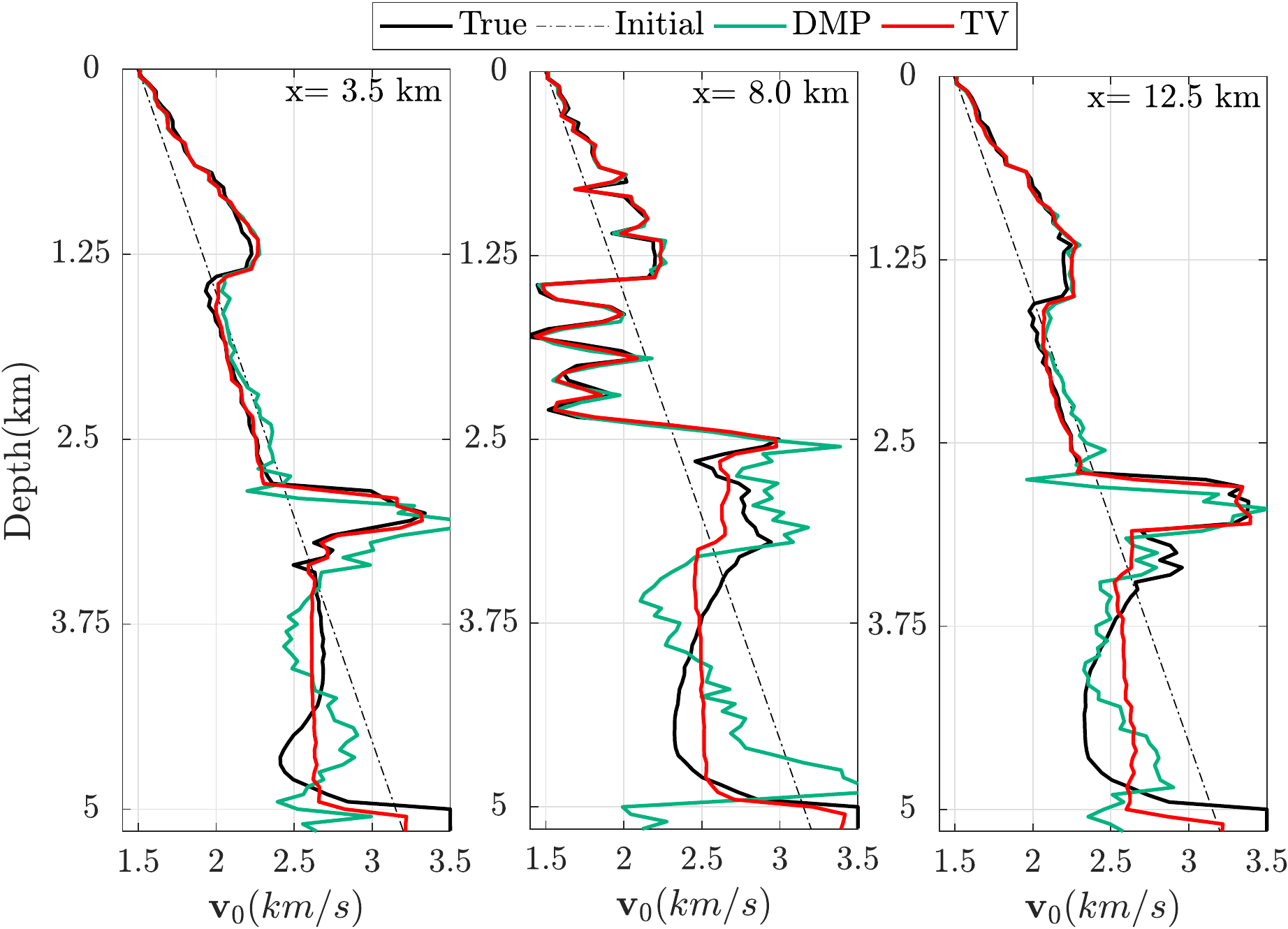} 
\caption{North Sea case study with noiseless data. Direct comparison along the logs at $x=3.5$ (left), $x=8.0$ (center) and $x=12.5 ~km$ (right) between the true velocity model (black), the initial model (dashed line) and the bound constrained IR-WRI models shown in Fig. \ref{fig:val_mono_noiseless}. The logs of the IR-WRI models obtained with DMP regularization and DMP+TV regularization are the green lines and the red lines, respectively.}
\label{fig:val_mono_logs_noiseless}
\end{figure}

We start with mono-parameter bound-constrained IR-WRI when $\bold{v}_0$ is the optimization parameter and $\bold{\delta}$ and $\bold{\epsilon}$ are passive parameters. We first consider noiseless data. \\
The final $\bold{v}_0$ models inferred from bound-constrained IR-WRI with DMP regularization only and DMP+TV regularizations are shown in Fig. \ref{fig:val_mono_noiseless}. Also, direct comparisons between the logs extracted from the true model, the initial model, and the IR-WRI models at $x=3.5~km$, $x=8.0~km$ and $x=12.5~km$ are shown in Fig.~\ref{fig:val_mono_logs_noiseless}.   
Although the crude initial $\bold{v}_0$ model and the smooth $\delta$ and $\epsilon$ passive models, the shallow sedimentary part and the gas layers are fairly well reconstructed with the two regularization settings. The main differences are shown at the reservoir level and below. Without TV regularization, the reconstruction at the reservoir level is quite noisy and the inversion fails to reconstruct the smoothly-decreasing velocity below the reservoir due to the lack of diving wave illumination at these depths. This in turn prevents the focusing of the deep reflector at 5~km depth by migration of the associated short-spread reflections.
When TV regularization is used, IR-WRI provides a more accurate and cleaner image of the reservoir and better reconstructs the sharp contrast on top of it. It also reconstructs the deep reflector at the correct depth in the central part of the model, while the TV regularization has replaced the smoothly-decreasing velocities below the reservoir by a piecewise constant layer with a mean velocity (Fig. \ref{fig:val_mono_logs_noiseless}).

%
%
%
%
\begin{figure}[ht!]
\centering
   \includegraphics[width=0.46\textwidth]{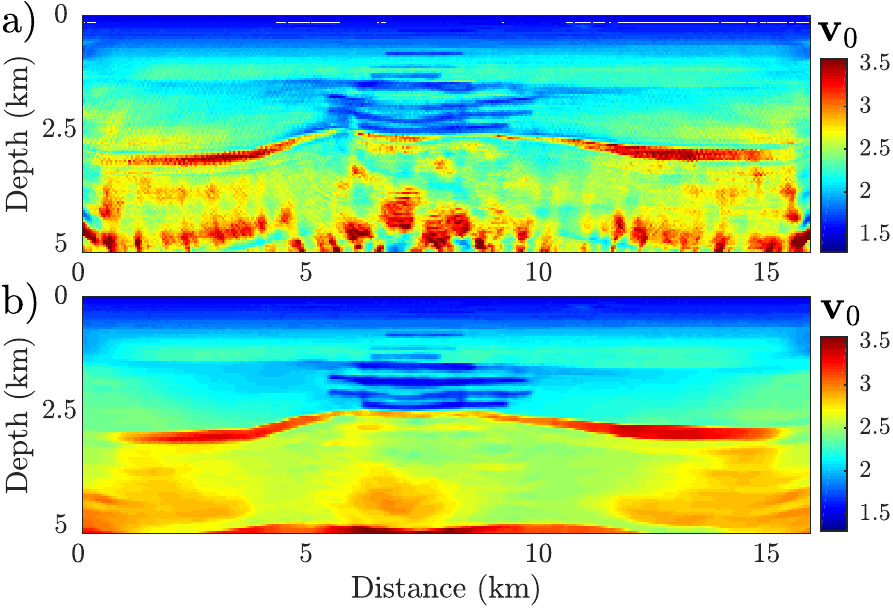} 
\caption{North Sea case study with noisy data (SNR=10~db). Mono-parameter bound constrained IR-WRI with (a) DMP (b) TV regularization.}
\label{fig:val_mono_noisy}
\end{figure}

\begin{figure}[ht!]
   \includegraphics[width=0.48\textwidth]{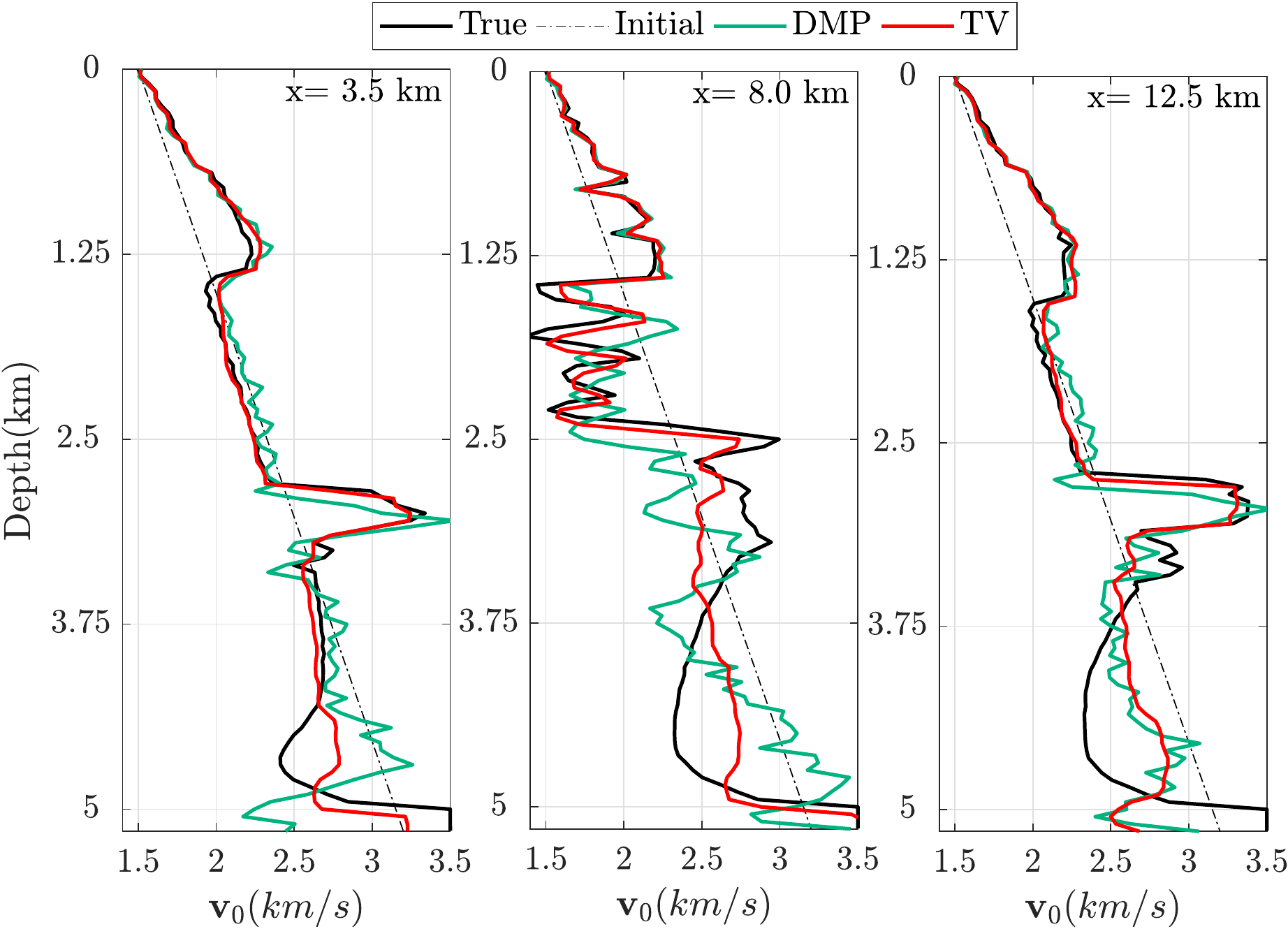} 
\caption{North Sea case study with noisy data (SNR=10~db). Direct comparison along the logs at $x=3.5$ (left), $x=8.0$ (center) and $x=12.5 ~km$ (right) between the true velocity model (black), the initial model (dashed line) and the bound constrained IR-WRI estimated models with DMP, Fig. \ref{fig:val_mono_noiseless}a, (green) and TV regularization, Fig. \ref{fig:val_mono_noiseless}b, (red).}
\label{fig:val_mono_logs_noisy}
\end{figure}

%
%
To generate more realistic test, we add random noise with Gaussian distribution to the data with SNR=10~db and re-run the mono-parameter bound-constrained IR-WRI with DMP and TV regularizations (Fig. \ref{fig:val_mono_noisy}). The direct comparison between the true model, the initial model and the IR-WRI models of Fig. \ref{fig:val_mono_noisy} at distances $x=3.5~km$, $x=8.0~km$ and $x=12.5~km$ are shown in Fig. \ref{fig:val_mono_logs_noisy}. The noise degrades the reconstruction of the gas layers both in terms of velocity amplitudes and positioning in depth when only DMP regularization is used. Also, the sharp reflector on top of the reservoir is now unfocused and mis-positioned in depth accordingly (Fig.~\ref{fig:val_mono_logs_noisy}). The TV regularization significantly reduces these amplitude and mis-positioning errors (Figs. \ref{fig:val_mono_noisy}b and \ref{fig:val_mono_logs_noisy}) and hence produces $\bold{v}_0$ models which are much more consistent with those obtained for noiseless data (Figs. \ref{fig:val_mono_noiseless}b and \ref{fig:val_mono_logs_noiseless}).

\begin{figure}[ht!]
\centering
   \includegraphics[width=0.48\textwidth]{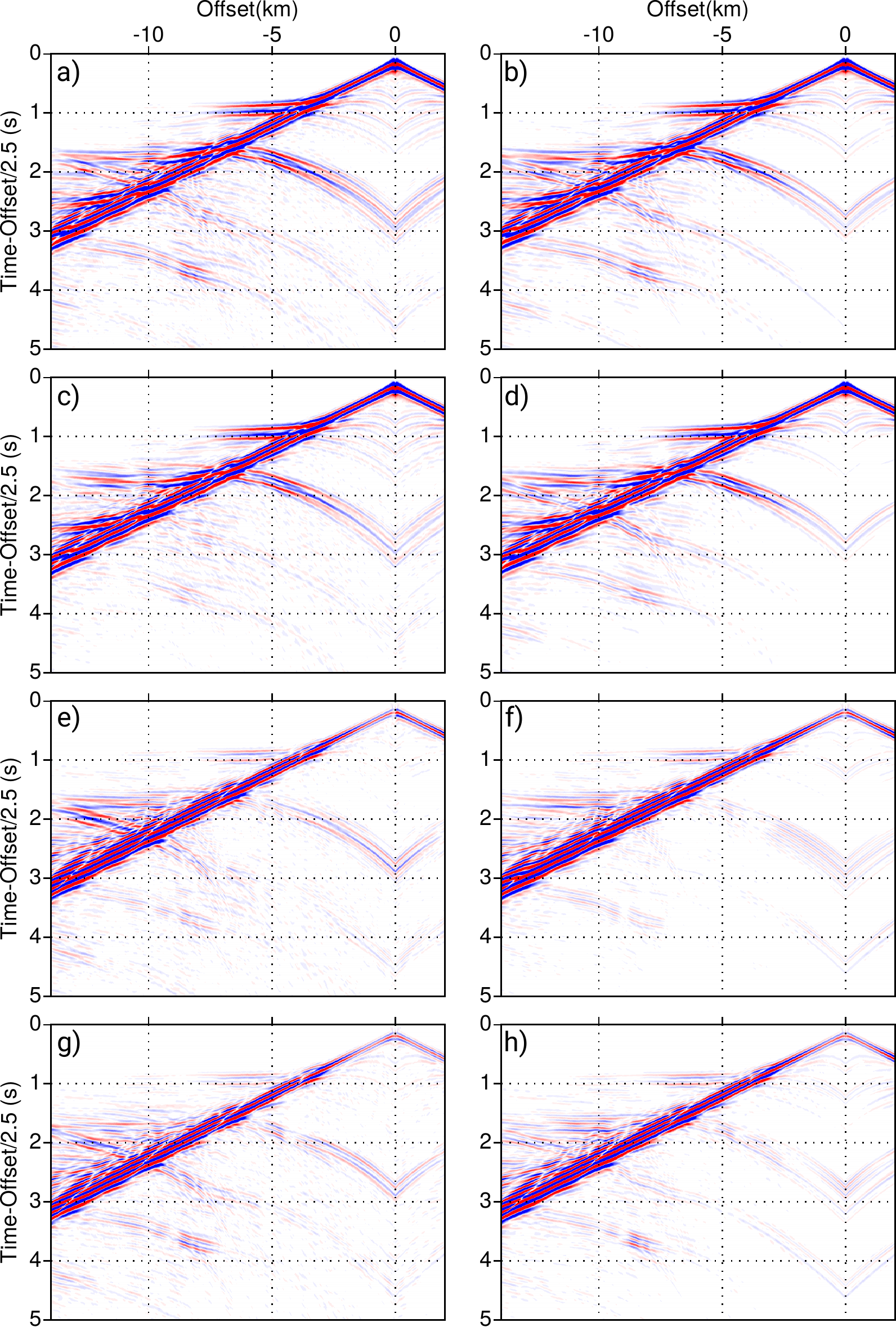} 
\caption{North Sea case study. Time-domain seismograms computed in the subsurface models obtained by ($\bold{v}_0$) mono-parameter IR-WRI. (a-b) Noiseless data: (a) DMP regularization, (b) DMP+TV regularization. (c-d) Same as (a-b) for noisy data.
(e-h) Residuals between the seismograms computed in the true model (Fig. \ref{fig:val_seismo_true}a) and those shown in (a-d). The seismograms are plotted with a reduction velocity of 2.5 km/s.}
\label{fig:sismos_mono}
\end{figure}

\begin{figure*}[ht!]
\centering
   \includegraphics[width=0.92\textwidth]{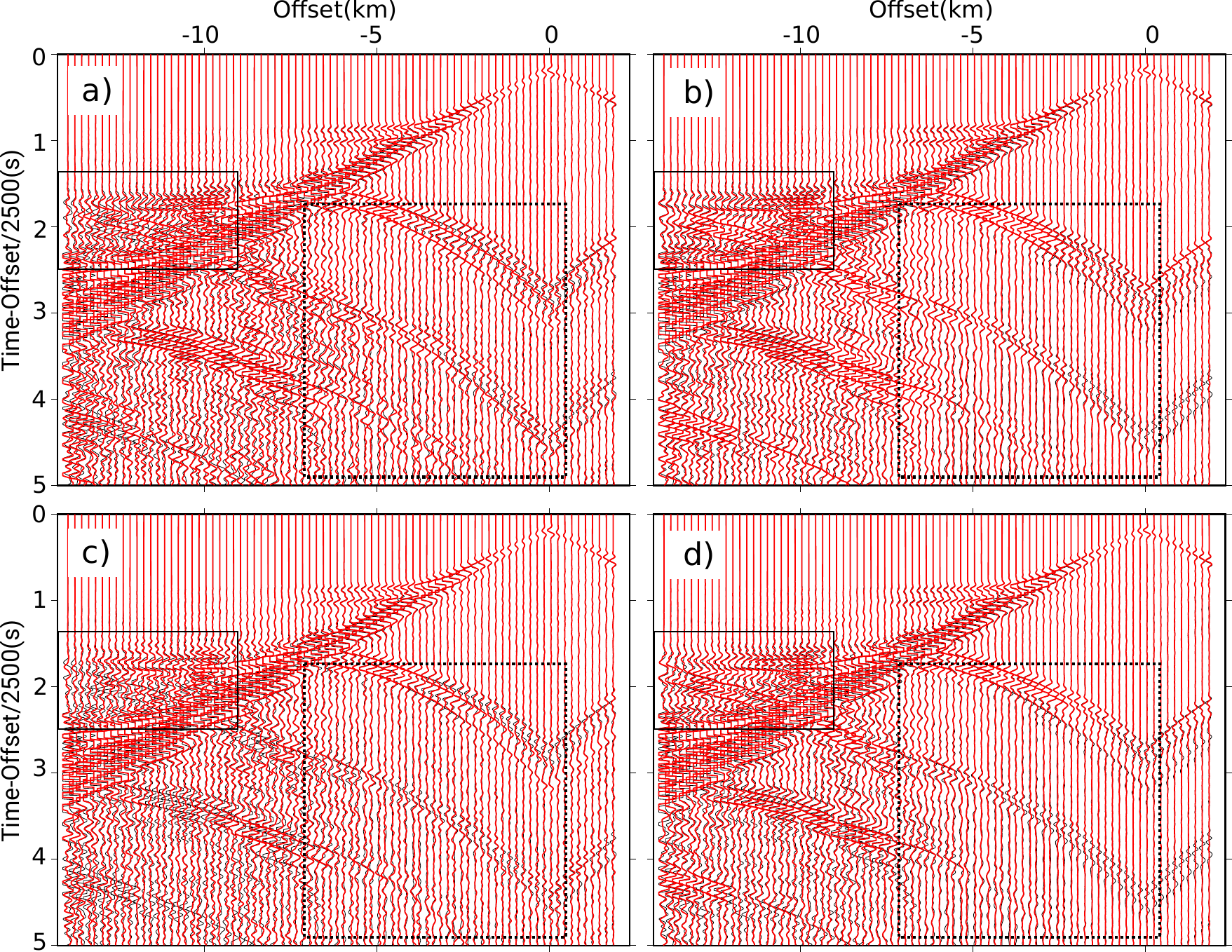}
\caption{North Sea case study. Mono-parameter bound-constrained IR-WRI results. Direct comparison between seismograms computed in the true model (black) and in the IR-WRI models (red). (a-b) Noiseless data (IR-WRI models of Fig. \ref{fig:val_mono_noiseless}). (a) DMP regularization. (b) DMP+TV regularization. (c-d) Same as (a-b) for noisy data (IR-WRI models of Fig.  \ref{fig:val_mono_noisy}). The seismograms are plotted with a reduction velocity of 2.5 km/s. True amplitudes are shown after a gain with offset and time for amplitude balancing. The solid box delineates post-critical reflections from the reservoir and the refracted wave from the deep interface, while the dot boxes delineate pre-critical reflections from the reservoir and the deep reflector. The amplitudes of the reverberating guided waves at long offsets are clipped for sake of clarity.}
\label{fig:sismos_mono_direct}
\end{figure*}

To assess how the differences between the velocity models shown in Fig. \ref{fig:val_mono_noiseless} impact waveform match, we compute time-domain seismograms in these models as well as the differences with those computed in the true model (Figs.~\ref{fig:sismos_mono} and \ref{fig:sismos_mono_direct}). The time-domain seismograms and the residuals shown in Fig.~\ref{fig:sismos_mono} give an overall vision of the achieved data fit, while the direct comparison between the recorded and modelled seismograms shown in Fig. \ref{fig:sismos_mono_direct} allow for a more detailed assessment of the waveform match for a specific arrival. 

A first conclusion is that, for all of the models shown in Fig. \ref{fig:val_mono_noiseless}, the main arrivals, namely those which have a leading role in the reconstruction of the subsurface model (diving waves, pre- and post-critical reflections), are not cycle skipped relative to those computed in the true model. We note however more significant residuals for the reverberating guided waves when TV regularization is used. These mismatches were generated by small wavespeed errors generated in the shallow part of the model by the TV regularization. These artifacts can be probably corrected by deactivating or decreasing the weight of the TV regularization locally. 
For noiseless data, the direct comparison between the seismograms computed in the true model and in the reconstructed ones show how the DMP+TV regularization improves the waveform match both at pre- and post-critical incidences relative to the DMP regularization alone (Fig. \ref{fig:sismos_mono_direct}a,b). 
For noisy data, the data fit is slightly degraded by noise when DMP regularization is used, while the TV regularization produces a data fit which is more consistent with that obtained with noiseless data (Compare Figs. \ref{fig:sismos_mono_direct}a and \ref{fig:sismos_mono_direct}c for DMP regularization, and Figs. \ref{fig:sismos_mono_direct}b and \ref{fig:sismos_mono_direct}d for DMP+TV regularization). It is striking to see the strong impact of noise on the quality of the $\bold{v}_0$ reconstruction when DMP regularization is used (compare Figs. \ref{fig:val_mono_noiseless}a and \ref{fig:val_mono_noisy}a, and Figs. \ref{fig:val_mono_logs_noiseless} and \ref{fig:val_mono_logs_noisy}, green curves), compared to its more moderate impact on the data fit (Compare Figs. \ref{fig:sismos_mono_direct}a and \ref{fig:sismos_mono_direct}c). This highlights the ill-posedness of the FWI, which is nicely mitigated by the prior injected by the TV regularization as illustrated by the consistency of the $\bold{v}_0$ models inferred from noiseless and noisy data (Compare Figs. \ref{fig:sismos_mono_direct}b and \ref{fig:sismos_mono_direct}b).
%
%
%
%
%
%
%

%
%
%
%
%
%
\subsubsection{Multi-parameter IR-WRI}
\begin{figure}[ht!]
\centering
      \includegraphics[width=0.48\textwidth]{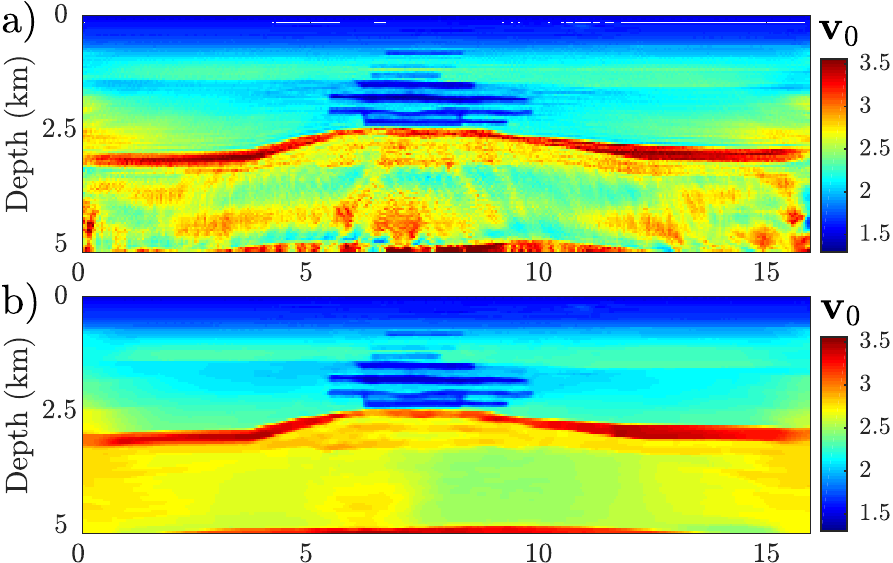} \\
   \includegraphics[width=0.487\textwidth]{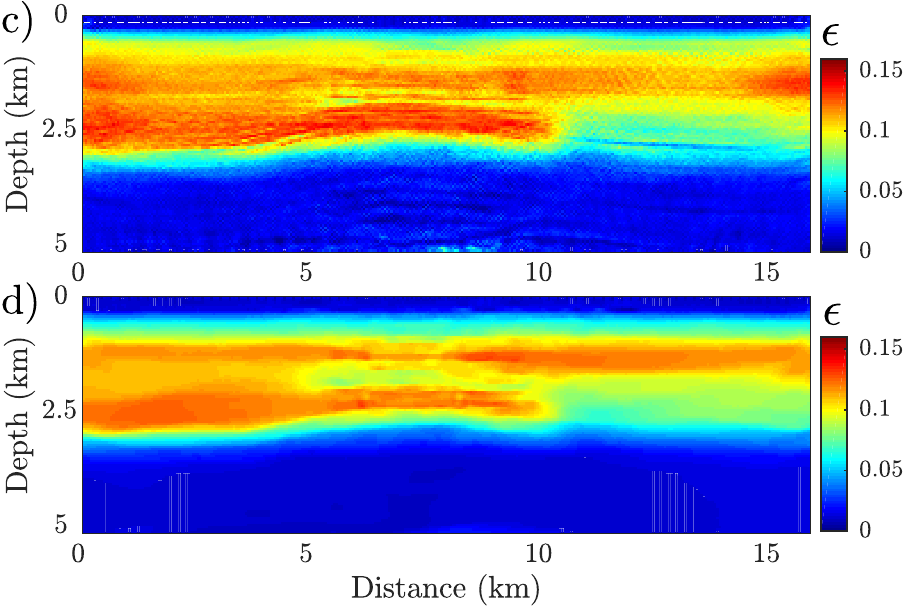}
\caption{North Sea case study with noiseless data. Multi-parameter IR-WRI results. (a-b) $\bold{v}_0$ models obtained with (a) DMP (a) and (b) DMP+TV regularizations. (c-d) Same as (a-b) for $\epsilon$.}
\label{fig:val_joint_noiseless}
\end{figure}

\begin{figure}[ht!]
   \includegraphics[width=0.48\textwidth]{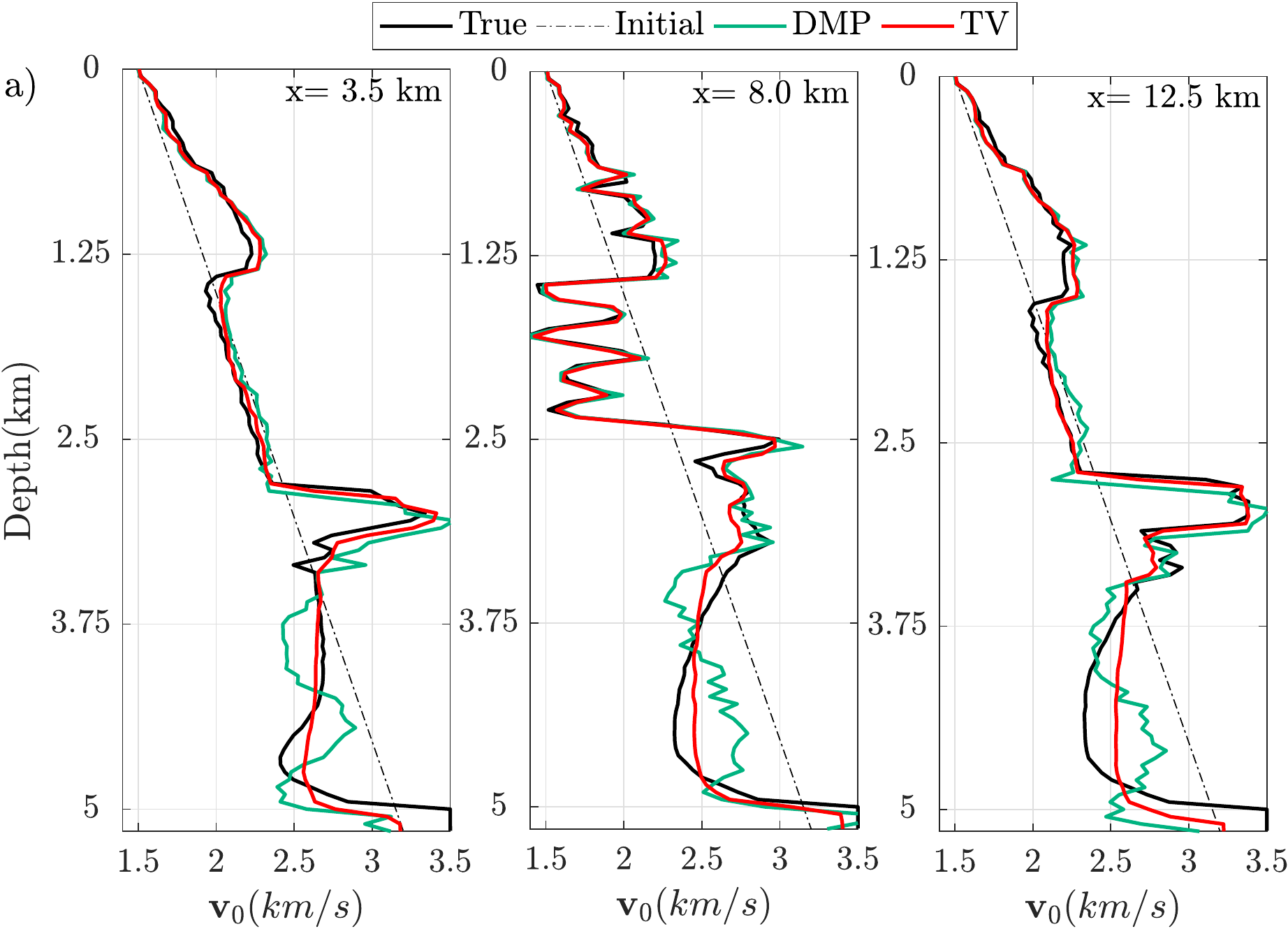} \\
   \includegraphics[width=0.47\textwidth]{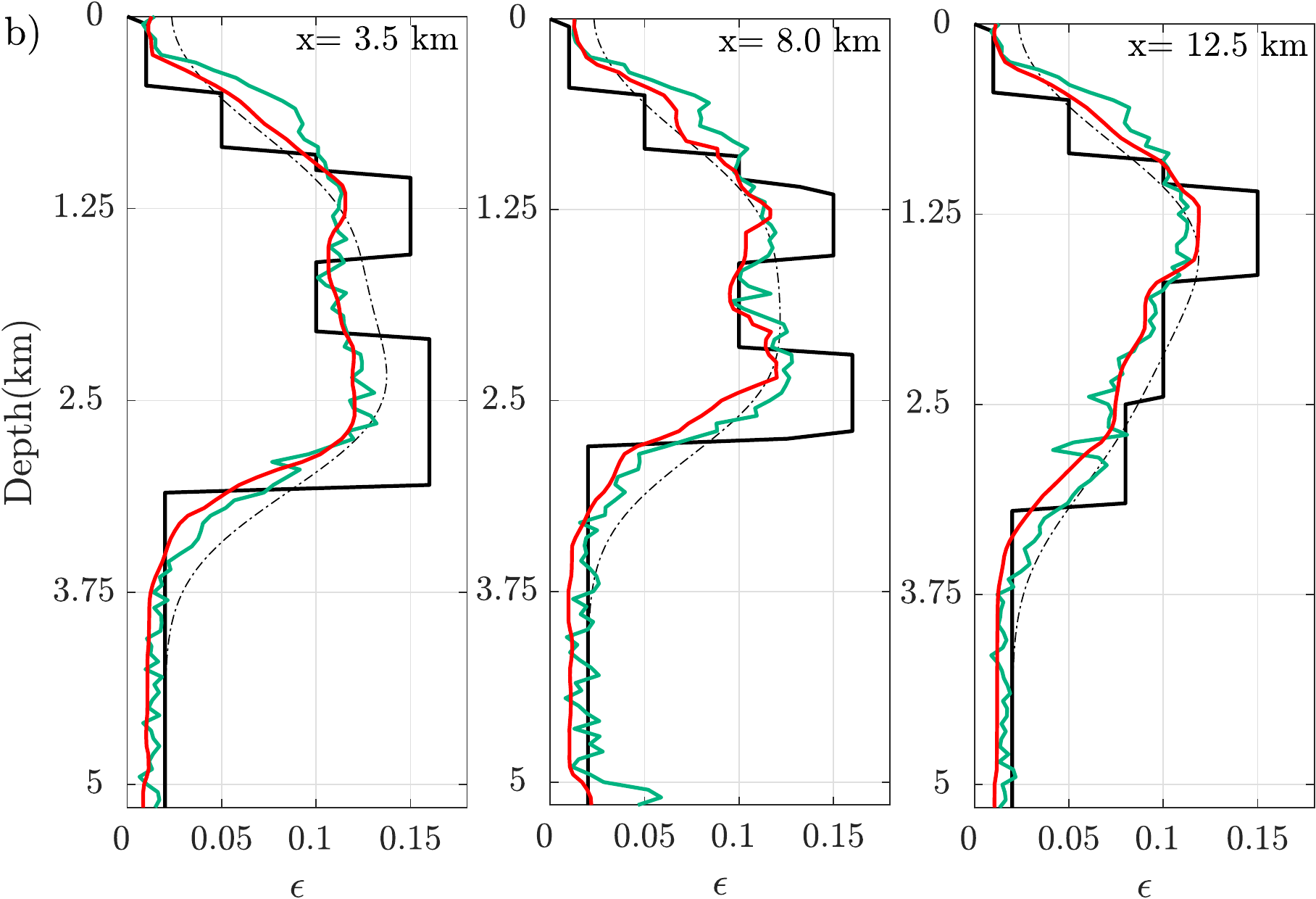}
\caption{North Sea case study with noiseless data. Multi-parameter bound-constrained IR-WRI results.  Direct comparison  along logs at $x=3.5$ (left), $x=8.0$ (center) and $x=12.5 ~km$ (right) between the true model (black), the initial model (dashed line) and the models inferred from IR-WRI with DMP (green) and DMP+TV regularization (red) for (a) $\bold{v}_0$ and (b) $\bold{\epsilon}$ (Fig. \ref{fig:val_joint_noiseless}).}
\label{fig:val_joint_logs_noiseless}
\end{figure}

We continue with multi-parameter bound-constrained IR-WRI when $\bold{v}_0$ and $\bold{\epsilon}$ are optimization parameters and $\bold{\delta}$ is a passive parameter. As for the mono-parameter inversion, we start with noiseless data.
The final $\bold{v}_0$ and $\bold{\epsilon}$ models inferred from bound-constrained IR-WRI with DMP and DMP+TV regularizations are shown in Fig. \ref{fig:val_joint_noiseless}. The direct comparisons between the logs extracted from the true $\bold{v}_0$ model, the initial model, the bound constrained IR-WRI models with DMP and DMP+TV regularization at $x=3.5~km$, $x=8.0~km$ and $x=12.5~km$ are shown in Fig. \ref{fig:val_joint_logs_noiseless}a, and the same comparisons for $\bold{\epsilon}$ are depicted in Fig. \ref{fig:val_joint_logs_noiseless}b. Compared to the mono-parameter inversion results, involving $\epsilon$ as an optimization parameter clearly improves the reconstruction at the reservoir level down to around 3~km depth (compare Figs. \ref{fig:val_mono_noiseless} and \ref{fig:val_joint_noiseless}a,b).
The long to intermediate wavelengths of the $\epsilon$ model are primarily updated according to the radiation pattern of this parameter in the ($\bold{v}_0,\epsilon,\delta$) parametrisation. The main effect of TV regularization relative to DMP regularization is to remove high-frequency noise from the $\epsilon$ model (Fig. \ref{fig:val_joint_logs_noiseless}b).

The time-domain seismograms computed in the multi-parameter models inferred from bound-constrained IR-WRI when DMP and DMP+TV regularizations are used are shown in Fig. \ref{fig:val_seismo_joint}. The direct comparison between the recorded and modelled seismograms is shown in Fig. \ref{fig:sismos_joint_direct}a,b.
Clearly, using both $\bold{v}_0$ and $\epsilon$ as optimization parameters allows us to better conciliate the fit of the pre- and post-critical reflections (Compare Figs. \ref{fig:sismos_mono_direct}a and \ref{fig:sismos_joint_direct}a for DMP regularization, and Figs. \ref{fig:sismos_mono_direct}b and \ref{fig:sismos_joint_direct}b for DMP+TV regularization). With noiseless data, DMP regularization is enough to achieve a high-fidelity data fit, which looks better than that obtained with DMP+TV regularization (Compare Figs.  \ref{fig:sismos_joint_direct}a and \ref{fig:sismos_joint_direct}b). 
%
%
%
%
%
%

\begin{figure}[ht!]
\centering
      \includegraphics[width=0.48\textwidth]{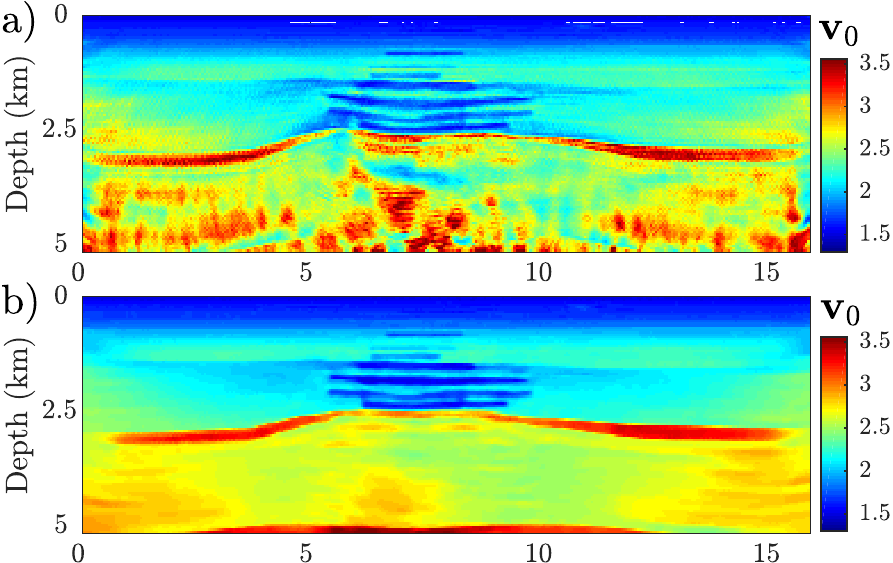} \\
   \includegraphics[width=0.487\textwidth]{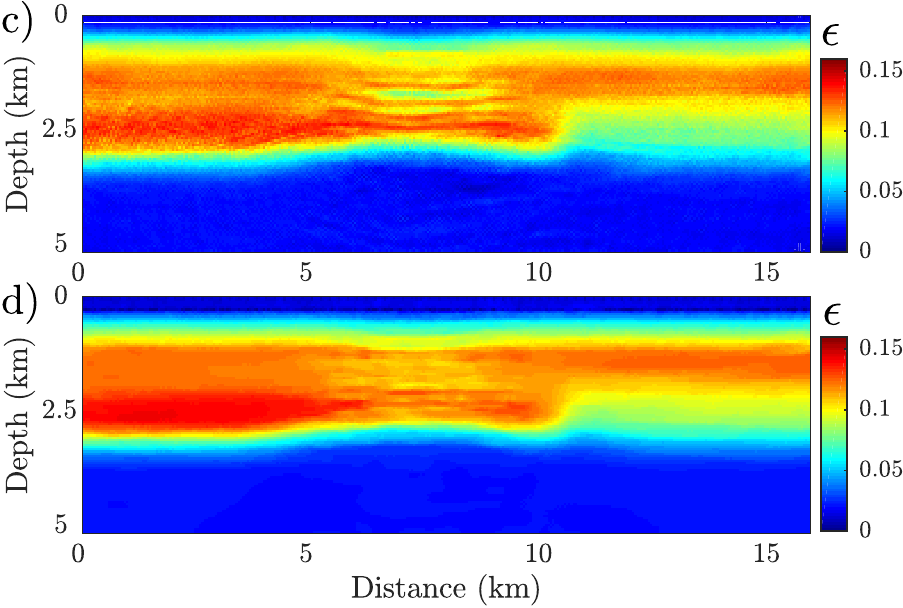}
\caption{North Sea case study with noisy data (SNR=10~db). Multi-parameter bound-constrained IR-WRI results. (a-b) $\bold{v}_0$ models obtained with (a) DMP (a) and (b) DMP+TV regularizations. (c-d) Same as (a-b) for $\epsilon$.}
\label{fig:val_joint_noisy}
\end{figure}

\begin{figure}[ht!]
   \includegraphics[width=0.48\textwidth]{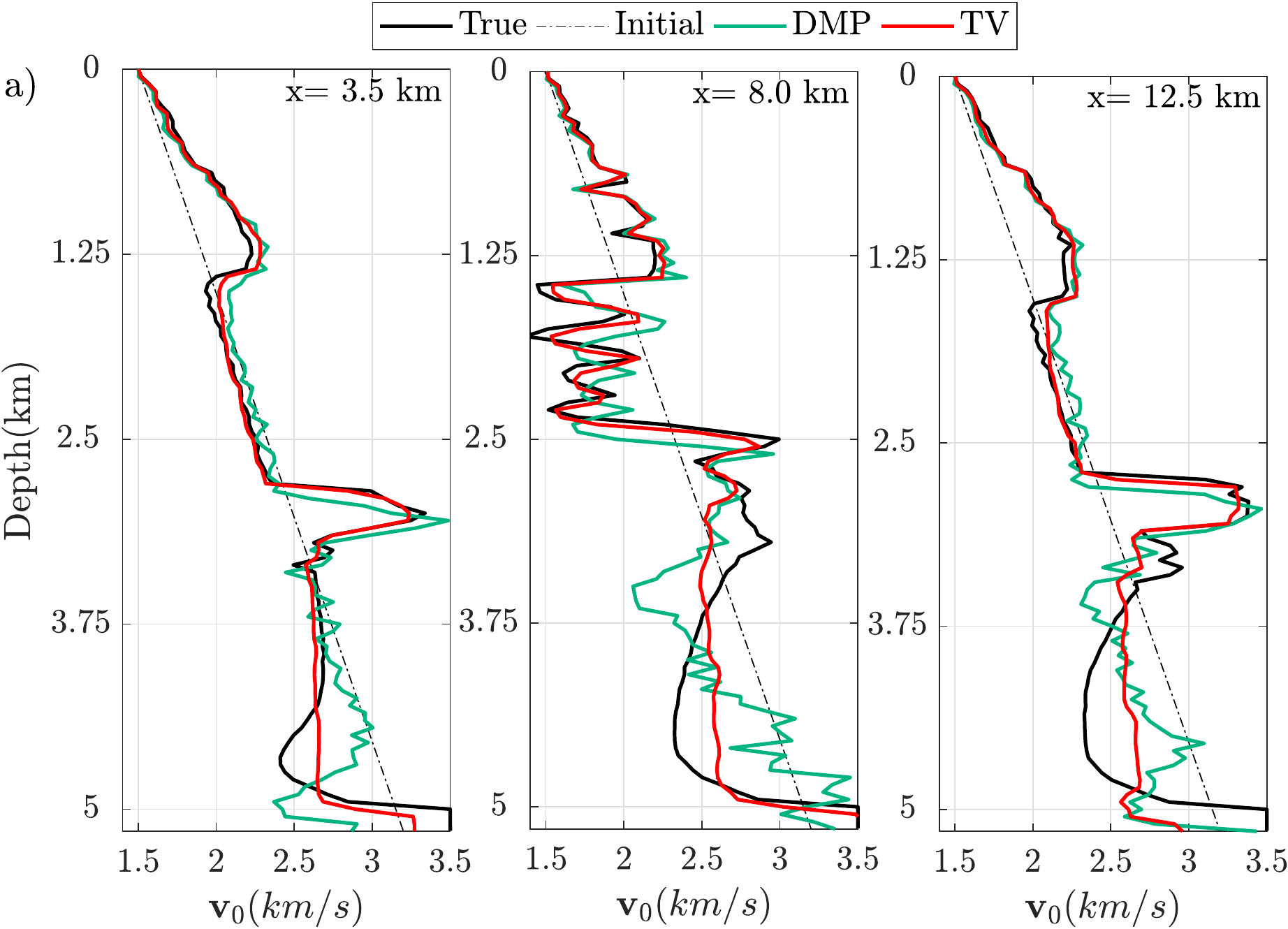} \\
   \includegraphics[width=0.47\textwidth]{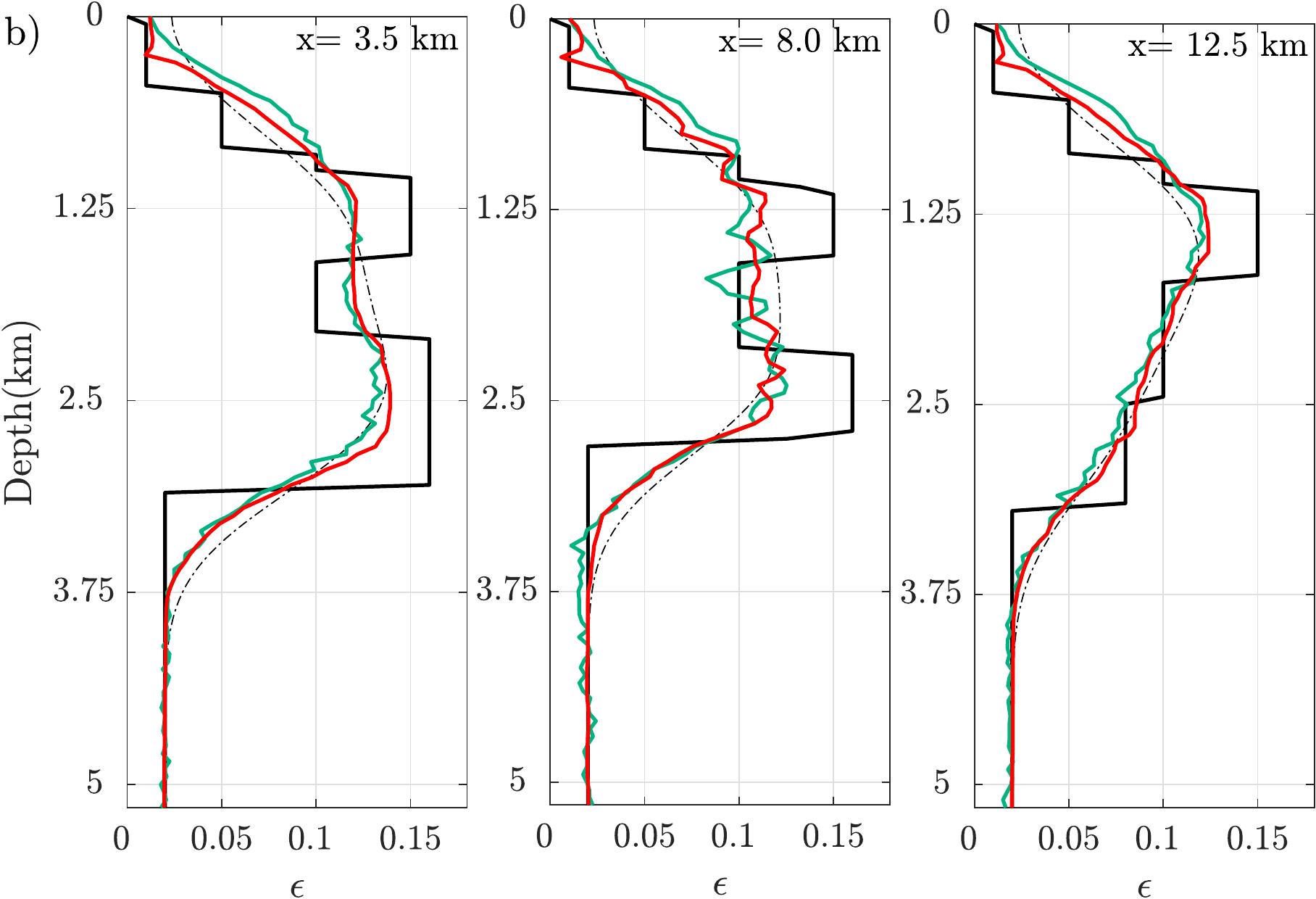}
\caption{North Sea case study with noisy data.  Multi-parameter bound-constrained IR-WRI results. Direct comparison  along logs at $x=3.5$ (left), $x=8.0$ (center) and $x=12.5 ~km$ (right) between the true model (black), the initial model (dashed line) and the models inferred from IR-WRI with DMP (green) and DMP+TV regularization (red) for (a) $\bold{v}_0$ and (b) $\bold{\epsilon}$ (Fig. \ref{fig:val_joint_noisy}).}
\label{fig:val_joint_logs_noisy}
\end{figure}

\begin{figure}[ht!]
\centering
   \includegraphics[width=0.48\textwidth]{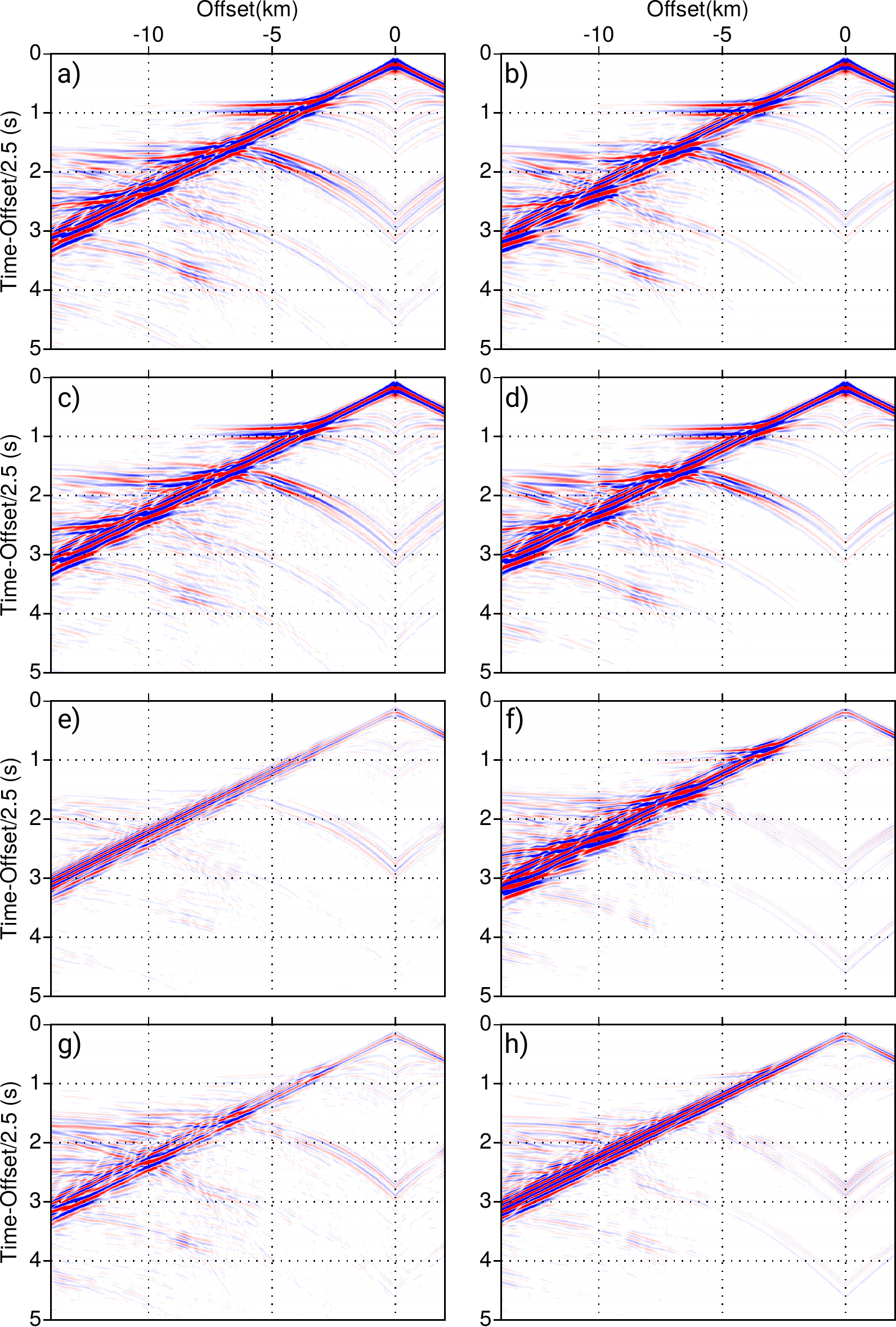} 
\caption{North Sea case study. Same as Fig. \ref{fig:sismos_mono} for ($v_0$,$\epsilon$) multi-parameter IR-WRI.}
\label{fig:val_seismo_joint}
\end{figure}

\begin{figure*}[ht!]
\centering
   \includegraphics[width=0.92\textwidth]{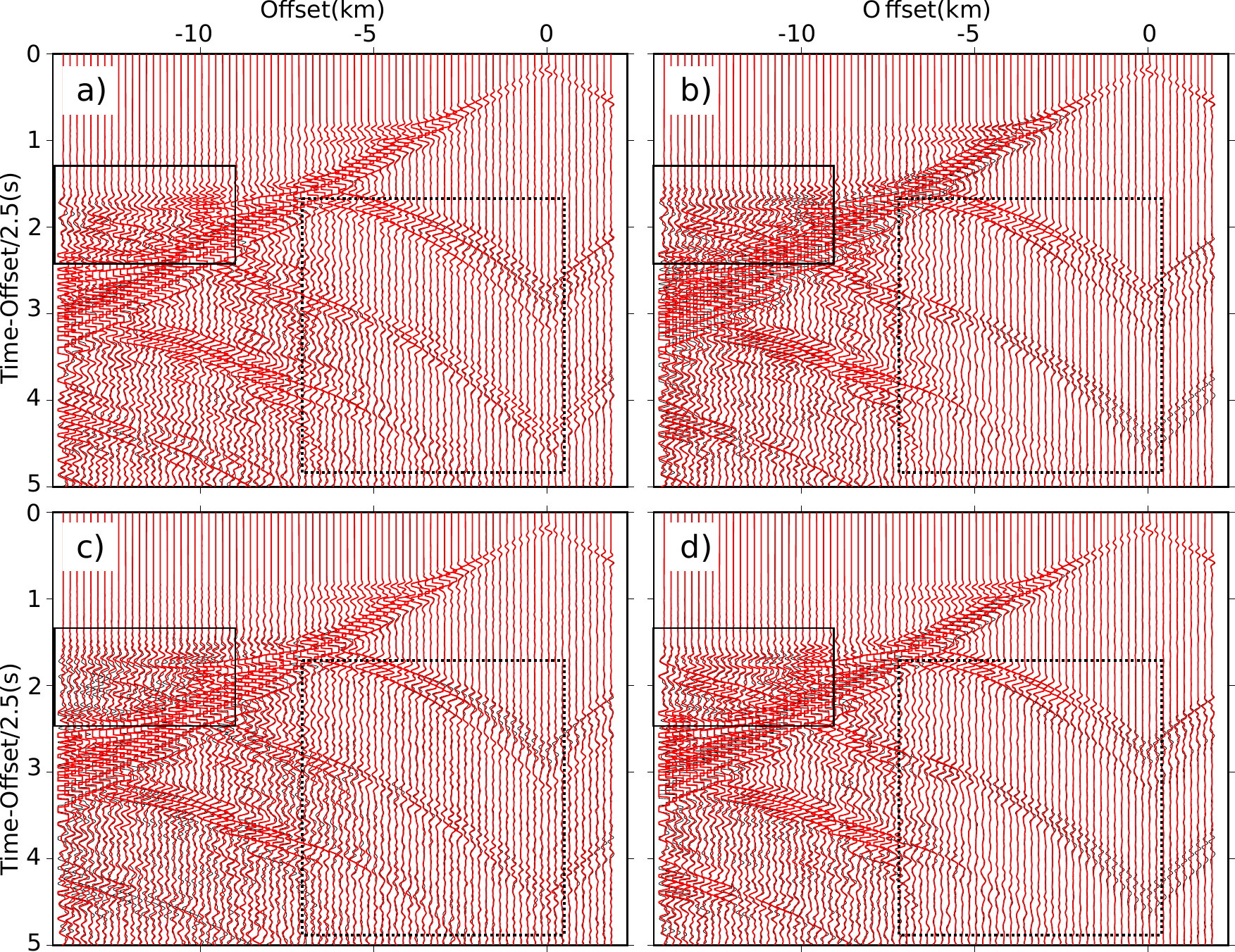} 
\caption{North Sea case study.  Same as Fig. \ref{fig:sismos_mono_direct} for ($v_0$,$\epsilon$) multi-parameter IR-WRI.}
\label{fig:sismos_joint_direct}
\end{figure*}

We repeat now the joint inversion when data are contaminated by Gaussian random noise with a SNR=10~db. The final $\bold{v}_0$ and $\bold{\epsilon}$ models inferred from bound-constrained IR-WRI with DMP and DMP+TV regularization are shown in Fig. \ref{fig:val_joint_noisy}.
Direct comparisons along logs extracted from the true models, the initial models, and the models inferred from the bound constrained multi-parameter IR-WRI with DMP and DMP+TV regularization at $x=3.5~km$, $x=8.0~km$ and $x=12.5~km$ are shown in Fig. \ref{fig:val_joint_logs_noisy}. For the $\bold{v}_0$ reconstruction, a trend similar to that shown for the mono-parameter inversion is shown, with a more significant impact of the noise on the velocity model reconstructed with DMP regularization compared to the one reconstructed with DMP+TV regularization (Compare Figs.~\ref{fig:val_joint_noiseless} and \ref{fig:val_joint_noisy}). As expected, the impact of noise is more significant on the second-order $\epsilon$ model, even when TV regularization is used, in the sense that the estimated perturbations of the initial model have smaller amplitudes compared to the noiseless case (Fig. ~ \ref{fig:val_joint_logs_noisy}b).

The time-domain seismograms are shown in Figs.~\ref{fig:val_seismo_joint} and \ref{fig:sismos_joint_direct}c,d. The data fit obtained with DMP regularization has been significantly degraded compared to that obtained with noiseless data (compare Figs. \ref{fig:sismos_joint_direct}a and \ref{fig:sismos_joint_direct}c), while the data fit obtained with DMP+TV regularization is more consistent with noiseless and noisy data (compare Figs. \ref{fig:sismos_joint_direct}b and \ref{fig:sismos_joint_direct}d). This is consistent with the previous conclusions drawn from the mono-parameter inversion.
 
%
%
%
%
%
%
%
\section{Discussion}
We have extended the ADMM-based wavefield reconstruction inversion method (IR-WRI), originally developed for mono-parameter wavespeed reconstruction \citep{Aghamiry_2019_IWR,Aghamiry_2019_IBC}, to multi-parameter inversion in VTI acoustic media. 
We have first discussed which formulations of the VTI acoustic wave equation are bilinear in wavefield and subsurface parameters. First-order velocity stress form is often more convenient than the second-order counterpart to fullfil bilinearity, in particular if density (or buoyancy) is an optimization parameter. However, it may be not the most convenient one for frequency-domain wavefield reconstruction as the size of the linear system to be solved scales to the number of wavefield components. To bypass this issue, wavefield reconstruction and parameter estimation can be performed with different wave equations, provided they give consistent solutions (\citet{Gholami_2013_WPA2} and this study). 

Bilinearity allows us to recast the parameter estimation subproblem as a linear subproblem and hence the waveform inversion as a biconvex problem. Although ADMM has been originally developed to solve distributed convex problems, it can be used as is to solve biconvex problem \citep{Boyd_2011_DOS}. From the mathematical viewpoint, biconvex problems should have superior convergence properties compared to fully nonconvex problem \citep[this discussion is out of the scope of this study but we refer the interested reader to][]{Benning_2015_ADM}. Alternatively, ADMM-like optimization can be used to perform IR-WRI heuristically without bilinear wave equation, hence keeping the parameter estimation subproblem, equation~\ref{main4_m}, nonlinear. Accordingly, equation~\ref{main4_m} would be solved with a Newton algorithm rather than with a Gauss-Newton one. 
This nonlinear updating of the parameters may however require several inner Newton iterations per IR-WRI cycle, while \citet{Aghamiry_2019_IWR} showed that one inner Gauss-Newton iteration without any line search was providing the most efficient convergence of IR-WRI when bilinearity is fulfilled. 

Indeed, the bilinearity specification limits the choice of subsurface parametrisation for parameter estimation. In the general case of triclinic elastodynamic equations, a subsurface parametrisation involving buoyancy and stiffness or compliance coefficients will be the most natural ones, as they correspond to the coefficients of the equation of motion and the Hooke's law.
In the particular case of the VTI acoustic wave equation, we have developed a bilinear wave equation whose coefficients depend on the vertival wavespeed $\bold{v}_0$ and the Thomsen's parameters $\delta$ and $\epsilon$. Although $\bold{v}_0$ and $\epsilon$ are coupled at wide scattering angles, the $(\bold{v}_0,\epsilon,\delta)$ parametrisation was promoted by \citet{Gholami_2013_WPA1} and \citet{Gholami_2013_WPA2} because the dominant parameter $\bold{v}_0$ has a radiation pattern which doesn't depend on the scatteting angle, and hence can be reconstructed with a high resolution from wide-azimuth long-offset data. The counterpart is that updating the long wavelengths of the secondary parameter $\epsilon$ is challenging and requires so far a crude initial guess of its long wavelengths (Fig.~\ref{fig:val_initial}e,f), which can be used as prior to regularize the $\epsilon$ update.
Comparing the results of IR-WRI when $\epsilon$ is used as a passive parameter and as an optimization parameter shows that the sensitivity of the inversion to $\epsilon$ remains small provided that a reasonable guess of its long wavelengths are provided in the starting model (for the models, compare Figs. \ref{fig:val_mono_noiseless}-\ref{fig:val_mono_noisy} and Figs. \ref{fig:val_joint_noiseless}-\ref{fig:val_joint_noisy}; for the data fit, compare Fig. \ref{fig:sismos_mono_direct} and Fig. \ref{fig:sismos_joint_direct}). This limited sensitivity of FWI to the short-to-intermediate wavelengths of $\epsilon$ in the $(\bold{v}_0,\bold{\epsilon},\bold{\delta})$ prompted for example \citet{Debens_2015_GAF} to estimate a crude $\epsilon$ background model with a coarse parametrisation by global optimization. 

Among the alternative parametrisations proposed for VTI acoustic FWI, \citet{Plessix_2011_GJI} proposed the $(\bold{v}_n,\bold{v}_h,\delta)$ or the $(\bold{v}_n,\bold{\eta},\delta)$ parametrisations for long-offset acquisition,  while \citet{Alkhalifah_2014_ANI} promoted the $(\bold{v}_h,\eta,\epsilon)$ parametrisation, where $\bold{v}_n = \bold{v}_0 \sqrt{1 + 2 \delta}$ is the so-called NMO velocity, $\bold{v}_h = \bold{v}_0 \sqrt{1 + 2 \epsilon}$ is the horizontal velocity and $\eta = (\epsilon-\delta)/(1 + 2 \delta)$ represents the anellipticity of the anisotropy. 
For the parametrisation promoted by  \citet{Plessix_2011_GJI}, the VTI equation developed by \citet{Zhou_2006_AAWa} given by
\begin{eqnarray}
\omega^2 \text{diag}\left(\frac{1}{\bold{v}_n^2}\right) \bold{u}_q + 2 \text{diag}(\bold{\eta}) \nabla_{\!\! xx} (\bold{u}_p + \bold{u}_q) &=& \bold{s}_q, \nonumber \\
 \omega^2  \text{diag}\left(\frac{1}{\bold{v}_n^2}\right) \bold{u}_p + \nabla_{\!\! xx} (\bold{u}_p + \bold{u}_q)  &+& \nonumber \\
   \text{diag}\left(\frac{1}{\sqrt{1+2 \bold{\delta}}}\right) \nabla_{\!\! zz} \text{diag}\left(\frac{1}{\sqrt{1+2 \bold{\delta}}}\right) \bold{u}_p &=& \bold{s}_p, 
\end{eqnarray}
is bilinear in wavefields and parameters $(1/\bold{v}_n^2,\eta)$, where $\bold{u}_p = \sqrt{1 + 2 \delta} \bold{u}_z $ and $\bold{u}_q = \bold{u}_x - \sqrt{1 + 2 \delta} \bold{u}_z $. Note that if $\delta$ is assumed to be smooth, the above equation can be approximated as
\begin{eqnarray}
\omega^2 \text{diag}\left(\frac{1}{\bold{v}_n^2}\right) \bold{u}_q + 2 \text{diag}(\bold{\eta}) \nabla_{\!\! xx} (\bold{u}_p + \bold{u}_q) &=& \bold{s}_q, \nonumber \\
\omega^2  \text{diag}\left(\frac{1}{\bold{v}_n^2}\right) \bold{u}_p + \nabla_{\!\! xx} (\bold{u}_p + \bold{u}_q)  &+& \nonumber \\
 \text{diag}\left(\frac{1}{1+2 \bold{\delta}}\right) \nabla_{\!\! zz} \bold{u}_p &=& \bold{s}_p 
\end{eqnarray}
which is bilinear in wavefields and parameters $(1/\bold{v}_n^2,\eta,1/(1+2\delta)$. This implies that $\delta$ can be involved as an optimization parameter if necessary. Note that, if this smoothness approximation is used to update the parameters during the primal problem, the modelled data and the source residuals can be solved with the exact equation to update the dual variables with more accuracy. 
For the $(\bold{v}_h,\eta,\epsilon)$ parametrisations \citep{Alkhalifah_2014_ANI}, according to \citet{Zhou_2006_AAWa} the VTI equations with smooth $\delta$ can be written as   
\begin{eqnarray}
\omega^2 \text{diag}\left(\frac{1}{\bold{v}_h^2}\right) \bold{u}_q + \nabla_{\!\! xx} (\bold{u}_p + \bold{u}_q) - \text{diag}\left(\frac{1}{1+2\bold{\eta}}\right) && \nonumber \\ \nabla_{\!\! xx} (\bold{u}_p + \bold{u}_q)&=& \bold{s}_q, \nonumber \\
\omega^2  \text{diag}\left(\frac{1}{\bold{v}_h^2}\right) \bold{u}_p + \text{diag}\left(\frac{1}{1+2\bold{\eta}}\right) \nabla_{\!\! xx} (\bold{u}_p + \bold{u}_q)  &+& \nonumber \\
 \text{diag}\left(\frac{1}{1+2 \bold{\epsilon}}\right) \nabla_{\!\! zz} \bold{u}_p &=& \bold{s}_p, \nonumber \\ 
\end{eqnarray}
which is bilinear in wavefields and parameters $(1/\bold{v}_h^2,1/(1+2\eta),1/(1+2\epsilon)).$
%
\section{Conclusion}
We have shown that ADMM-based IR-WRI can be extended to multi-parameter reconstruction for VTI acoustic media. The gradient of the misfit function for the parameter estimation subproblem involves the so-called virtual sources, which carry out the effect of the parameter-dependent radiation patterns. This suggests that, although IR-WRI extends the search-space of FWI to mitigate cycle skipping, it is impacted by ill-posedness associated with parameter cross-talk and incomplete angular illumination as classical FWI. We have verified this statement with a toy numerical example for which we have reproduced the same pathologies in terms of resolution and parameter cross-talks as those  produced by classical FWI during a former study. We have illustrated how equipping IR-WRI with bound constraints and TV regularization fully remove the ill-posedness effects  for this idealized numerical example. We have provided some guidelines to design bilinear wave equation of different order for different subsurface parametrisations. Although bilinearity puts some limitations on the choice of the subsurface parametrisation, it recasts the parameter estimation subproblem as a quadratic optimization problem, which can be solved efficiently with Gauss-Newton algorithm. 
Application on the long-offset synthetic case study representative of the North Sea has shown how IR-WRI can be started from a crude laterally-homogeneous vertical velocity model without impacting the inversion with cycle skipping, when a smooth $\delta$ parameter is used as a passive background model. However, a smooth initial $\epsilon$ model, albeit quite crude, is necessary to guarantee the convergence of the method to a good solution, either when $\epsilon$ is used as a passive or as an optimization parameter. For this case study where the low velocity gas layers and the smooth medium below the reservoir suffer from a deficit illumination of diving waves, the TV regularization plays a key role to mitigate the ill-posedness.

\section{Acknowledgments}
This study was partially funded by the SEISCOPE consortium (\textit{http://seiscope2.osug.fr}), sponsored by AKERBP, CGG, CHEVRON, EQUINOR, EXXON-MOBIL, JGI, PETROBRAS, SCHLUMBERGER, SHELL, SINOPEC and TOTAL. This study was granted access to the HPC resources of SIGAMM infrastructure (http://crimson.oca.eu), hosted by Observatoire de la C\^ote d'Azur and which is supported by the Provence-Alpes C\^ote d'Azur region, and the HPC resources of CINES/IDRIS/TGCC under the allocation 0596 made by GENCI.

%
\append{First and second-order wave equation with compliance notation} \label{Appc}
We consider the frequency-domain first-order velocity-stress equation in 2D VTI acoustic media with compliance notation as 
\begin{align} \label{edavti_comp}
-\hat{i} \omega \bold{v}_{x,l} &= \text{diag}(\bold{b})  \nabla_{\!\! x} \bold{u}_{x,l},  \nonumber \\
-\hat{i} \omega \bold{v}_{z,l} &= \text{diag}(\bold{b}) \nabla_{\!\! z}\bold{u}_{z,l},   \nonumber \\
-\hat{i} \omega [\text{diag}(\bold{\mathcal{s}}_{11}) \bold{u}_{x,l} +\text{diag}(\bold{\mathcal{s}}_{13}) \bold{u}_{z,l}]&= \nabla_{\!\! x}\bold{v}_{x,l}  - \hat{i} \omega \bold{s}_l \nonumber, \\
-\hat{i} \omega [\text{diag}(\bold{\mathcal{s}}_{13}) \bold{u}_{x,l}+\text{diag}(\bold{\mathcal{s}}_{33}) \bold{u}_{z,l}] &=  \nabla_{\!\! z} \bold{v}_{z,l} - \hat{i} \omega \bold{s}_l,
\end{align}
where $\bold{\mathcal{s}}_{ij} \in \mathbb{R}^{n \times 1}$ are the  compliance coefficients, and the other notations are defined after equation \ref{edavti}.
Gathering equation \ref{edavti_comp} for all sources results in the following matrix equation:
\begin{equation}
\label{eqvelstressf5_comp} 
 \begin{bmatrix}
\hat{i}\omega\bold{I} & \bold{B \nabla}\\
\bold{\nabla} & \hat{i}\omega\bold{\mathcal{S}}\\
\end{bmatrix} 
 \begin{bmatrix}
\bold{V}\\
\bold{U}\\
\end{bmatrix} 
 =   \hat{i}\omega
  \begin{bmatrix}
\bold{0}\\
\bold{S}\\
\end{bmatrix}, 
\vspace{-0.3cm}
\end{equation}
where 
\begin{equation}
\bold{\mathcal{S}}=
\begin{bmatrix}
\text{diag}(\bold{\mathcal{s}}_{11}) & \text{diag}(\bold{\mathcal{s}}_{13}) \\
\text{diag}(\bold{\mathcal{s}}_{13}) & \text{diag}(\bold{\mathcal{s}}_{33}) \\
\end{bmatrix}, \nonumber
\end{equation}
and the other notations are defined in equation \ref{eqvelstressf5}. 
Equation \ref{eqvelstressf5_comp} is linear in $\bold{U}$ and $\bold{V}$ when the model parameters embedded in $\bold{B}$ and $\bold{\mathcal{S}}$  are known. 
When $\bold{U}$ and $\bold{V}$  are known, this system can be recast as a new linear system in which the unknowns are the model parameters, and hence the bilinearity of the wave equation.
For the $l$th source, the new equation becomes
\begin{equation} \label{Eqval_m_comp}
\begin{bmatrix}
\bold{L}_{11} & \bold{0} & \bold{0} & \bold{0}   \\
\bold{L}_{21} & \bold{0} & \bold{0} & \bold{0} \\
 \bold{0} & \bold{L}_{32} & \bold{L}_{33}& \bold{0}  \\
 \bold{0} & \bold{0} & \bold{L}_{43}  & \bold{L}_{44}
\end{bmatrix}
\begin{bmatrix}
\bold{b} \\
\bold{\mathcal{s}}_{11} \\
\bold{\mathcal{s}}_{13} \\
\bold{\mathcal{s}}_{33} \\
\end{bmatrix}
=  
\begin{bmatrix}
\bold{0}-\hat{i} \omega\bold{v}_{x,l} \\
\bold{0}-\hat{i} \omega\bold{v}_{z,l} \\
 {\hat{i} \omega\bold{s}_l}-\nabla_{\!\! x}\bold{v}_{x,l}\\
{\hat{i} \omega\bold{s}_l}-\nabla_{\!\! z}\bold{v}_{z,l}\\
\end{bmatrix},
\end{equation}
where
\begin{equation}
\begin{cases}
\bold{L}_{11}=\text{diag}(\nabla_{\!\! x} \bold{u}_{x,l}), \\
 \bold{L}_{21}=\text{diag}(\nabla_{\!\! z} \bold{u}_{z,l}), \\
\bold{L}_{32}=\bold{L}_{43}=\hat{i} \omega\text{diag}(\bold{u}_{x,l}),\\
\bold{L}_{33}=\bold{L}_{44}=\hat{i} \omega\text{diag}(\bold{u}_{z,l}).
 \end{cases} \nonumber
\end{equation}
To develop the second-order wave equation, we eliminate $\bold{v}_{x,l}$ and $\bold{v}_{z,l}$ from equation \ref{edavti_comp}. We obtain the following equation
\begin{equation} \label{sovti_comp}
\begin{split}
&   \frac{1}{\omega^2}\nabla_{\!\! x} \text{diag}(\bold{b})  \nabla_{\!\! x}  \bold{u}_{x,l} + \text{diag}(\bold{\mathcal{s}}_{11}) \bold{u}_{x,l}+\text{diag}(\bold{\mathcal{s}}_{13}) \bold{u}_{z,l}= \bold{s}_l , \\
&  \frac{1}{\omega^2} \nabla_{\!\! z} \text{diag}(\bold{b})  \nabla_{\!\! z}  \bold{u}_{z,l} + \text{diag}(\bold{\mathcal{s}}_{13}) \bold{u}_{x,l}+ \text{diag}(\bold{\mathcal{s}}_{33}) \bold{u}_{z,l}= \bold{s}_l, \\
\end{split}
\end{equation}
which is bilinear with respect to buoyancy, compliance parameters and pressure wavefields.
With known buoyancy and compliance parameters, we get the following $2n \times 2n$ linear system to estimate wavefields for all sources
\begin{equation}
\label{linearu_comp}
\bold{A}(\bold{m}) \bold{U} = \bold{S},
\end{equation}
where $\bold{A}$ is given by
\begin{equation*}
\begin{bmatrix}
     \frac{1}{\omega^2}\nabla_{\!\! x} \text{diag}(\bold{b})  \nabla_{\!\! x}+\text{diag}(\bold{\mathcal{s}}_{11}) & \text{diag}(\bold{\mathcal{s}}_{13}) \\
    \text{diag}(\bold{\mathcal{s}}_{13})  & \frac{1}{\omega^2} \nabla_{\!\! z} \text{diag}(\bold{b})  \nabla_{\!\! z}+\text{diag}(\bold{\mathcal{s}}_{33})  
\end{bmatrix}.
\end{equation*}
and
\begin{equation}
\bold{m}=
 \begin{bmatrix}
\bold{b}\\
\bold{\mathcal{s}}_{11}\\
\bold{\mathcal{s}}_{13}\\
\bold{\mathcal{s}}_{33}\\
\end{bmatrix}, \nonumber
\end{equation}
and $\bold{S}$ is defined in equation \ref{UVS}.
When $\bold{U}$ is known, this system can also be recast as a new linear system in which the unknowns are the model parameters as 
\begin{equation}
\label{linearmaC}
\begin{bmatrix}
\bold{L}_{1} \\
\vdots \\
\bold{L}_{l} \\
\vdots \\
\bold{L}_{n_s}
\end{bmatrix}
\bold{m}
=
\begin{bmatrix}
\bold{y}_1 \\
\vdots \\
\bold{y}_l \\
\vdots \\
\bold{y}_{n_s}
\end{bmatrix},
\end{equation}
where $\bold{L}_l$ is given by 
\begin{equation*}
\begin{bmatrix}
    \frac{1}{\omega^2}\nabla_{\!\! x} \text{diag}(\nabla_{\!\! x} \bold{u}_{x,l})   & \text{diag}(\bold{u}_{x,l})  & \text{diag}( \bold{u}_{z,l})  & \bold{0}\\
  \frac{1}{\omega^2}\nabla_{\!\! z} \text{diag}(\nabla_{\!\! z} \bold{u}_{z,l})  & \bold{0} &  \text{diag}( \bold{u}_{x,l}) &    \text{diag}( \bold{u}_{z,l})
\end{bmatrix}
\end{equation*}
and
\begin{equation}
\bold{y}_{l}=
 \begin{bmatrix}
   \bold{s}_l  \\
    \bold{s}_l \\
\end{bmatrix}. \nonumber
\end{equation}

%
\append{Solving the optimization problem, equation \ref{main1_4}, with ADMM} \label{Appa}
Starting from an initial model $\bold{m}^0$, and setting the dual variables $\tilde{\bold{D}}$ and $\tilde{\bold{S}}$ equal to zero, the $k$th ADMM iteration for solving equation \ref{main1_4} reads as 
\citep[see][ for more details]{Aghamiry_2019_IWR,Aghamiry_2019_IBC,Aghamiry_2019_MAW}
\begin{subequations}
\label{main30}
\begin{eqnarray}
\bold{U}^{k+1} &\leftarrow &   \underset{\bold{U}}{\arg\min} ~~C_{\bold{m}^k,\tilde{\bold{D}}^k,\tilde{\bold{S}}^k}(\bold{U}), \label{main3a} \\
\bold{m}^{k+1} &\leftarrow &  \underset{\bold{m} \in \mathcal{C}}{\arg\min} ~~  C_{\bold{U}^{k+1},\tilde{\bold{S}}^k}(\bold{m}), \label{main3b} \\
\tilde{\bold{S}}^{k+1} &\leftarrow & \tilde{\bold{S}}^{k}  +\bold{S}- \bold{A}(\bold{m}^{k+1})\bold{U}^{k+1} , \label{main3c}\\
\tilde{\bold{D}}^{k+1} &\leftarrow & \tilde{\bold{D}}^{k} + \bold{D}- \bold{P} \bold{U}^{k+1} , \label{main3d}
\end{eqnarray} 
\end{subequations} 
where
\begin{equation} \label{Cu}
\begin{split}
C_{\bold{m}^k,\tilde{\bold{D}}^k,\tilde{\bold{S}}^k}(\bold{U})
&=  \lambda_0\|\bold{P}\bold{U}-\bold{D}-\tilde{\bold{D}}^k\|_{F}^2\\
&+\lambda_1\|\bold{A}(\bold{m}^k)\bold{U}-\bold{S}-\tilde{\bold{S}}^k\|_{F}^2 
 \end{split}
\end{equation}
and
\begin{equation} \label{Cm}
\begin{split}
C_{\bold{U}^{k+1},\tilde{\bold{S}}^k}(\bold{m})
& =\sum  \sqrt{|\partial_x \bold{m}|^2 + |\partial_z\bold{m}|^2}\\
&+\lambda_1\|\bold{A}(\bold{m})\bold{U}^{k+1}-\bold{S}-\tilde{\bold{S}}^k\|_{F}^2.
 \end{split}
\end{equation}
%
The regularized wavefields $\bold{U}$ are the minimizers of the quadratic cost function $C_{\bold{m}^k,\tilde{\bold{D}}^k,\tilde{\bold{S}}^k}(\bold{U})$, equation \ref{Cu}, where $\bold{m}^k$, $\tilde{\bold{D}}^k$ and $\tilde{\bold{S}}^k$ are kept fixed.
Zeroing the derivative of $C_{\bold{m}^k,\tilde{\bold{D}}^k,\tilde{\bold{S}}^k}(\bold{U})$ gives the wavefields as the solution of a linear system of equations defined by equation \ref{UE} (step 1 of the algorithm).
The regularized wavefields are then introduced as passive quantities in the cost function $C_{\bold{U}^{k+1},\tilde{\bold{S}}^k}(\bold{m})$, equation \ref{Cm}, which is minimized to estimate $\bold{m}$ over the desired set $\mathcal{C}$.
We solve this minimization subproblem with the splitting techniques. Accordingly, we introduce the auxiliary primal  variables $\bold{p}$ and $\bold{q}$ to decouple the $\ell{1}$ and $\ell{2}$ terms and split the parameter estimation subproblem into three sub-steps for $\bold{m}$, $\bold{p}$ and $\bold{q}$ \citep{Goldstein_2009_SBM,Aghamiry_2019_IBC}:
\begin{subequations}
\label{main4}
\begin{align}
\bold{m}^{k+1}&\leftarrow  \underset{\bold{m}}{\arg\min} ~~ 
\lambda_1\|\bold{A}(\bold{m})\bold{U}^{k+1}-\bold{S}-\tilde{\bold{S}}^k\|_{F}^2, \label{main4_m}\\
& +  \|\bold{\overline{\nabla} m}-\bold{p}^k-\tilde{\bold{p}}^k\|_{\bold{\Gamma}}^2 +
  \|\bold{m}-\bold{q}^k-\tilde{\bold{q}}^k\|_Z^2, 
  \nonumber\\
 \bold{p}^{k+1} &\leftarrow  \underset{\bold{p}}{\arg\min}  ~~\sum\sqrt{|\bold{p}_x|^2 + |\bold{p}_z|^2}  \label{main4_p}\\
 & \hspace{1.5cm}+  \|\overline{\nabla} \bold{m}^{k+1} -\bold{p}-\tilde{\bold{p}}^k\|_{\bold{\Gamma}}^2,  \nonumber\\
 \bold{q}^{k+1} &\leftarrow  \underset{\bold{q}\in \mathcal{C}}{\arg\min} ~~ \|\bold{m}^{k+1}-\bold{q}-\tilde{\bold{q}}^{k}\|_Z^2, \label{main4_q}
\end{align} 
\end{subequations}  
where $\overline{\nabla}$ is defined in equation \ref{partial},
\begin{equation}
\bold{p} = 
\begin{bmatrix}
\bold{p}_x\\
\bold{p}_z
\end{bmatrix},
\end{equation}
$\|\bold{x}\|^2_{\bold{\bullet}}=\bold{x}^T \bold{\bullet} \bold{x}$, $Z$ and $\bold{\Gamma}$ are diagonal matrices defined in equation \ref{Zeta} and \ref{Gamma}, respectively.

From the linearized equation \ref{linearma}, the subproblem for $\bold{m}$, equation \ref{main4_m}, can be written as
\begin{align} \label{quad}
\bold{m}^{k+1} &\leftarrow  \underset{\bold{m}}{\arg\min} ~~ 
\Biggr\Vert
\begin{bmatrix}
\bold{L}_{1} \\
\vdots\\
\bold{L}_{n_s} \\
\overline{\nabla} \\
\bold{I}
\end{bmatrix}
\bold{m} -
\begin{bmatrix}
\bold{y}^k_1+\tilde{\bold{s}}_1^k  \\
\vdots\\
\bold{y}^k_{n_s}+\tilde{\bold{s}}_{n_s}^k  \\  
\bold{p}^k+\tilde{\bold{p}}^k\\
\bold{q}^k+\tilde{\bold{q}}^k
  \end{bmatrix}
  \Biggr\lVert_{\bold{\Xi}}^2,
\end{align}
where the diagonal weighting matrix $\bold{\Xi}$ is defined as
 \begin{equation}
 \bold{\Xi}=
 \begin{bmatrix}
 \lambda_1\bold{I} & \bold{0} & \bold{0}\\
 \bold{0} & \Gamma& \bold{0}\\
  \bold{0} & \bold{0} & Z \\
\end{bmatrix} \in \mathbb{R}^{[3+n_s]n \times [3+n_s]n}. \nonumber
 \end{equation}
Equation \ref{quad} is now quadratic and admits a closed form solution as given in equation \ref{ME}  (step 2 of the algorithm).

The only remaining tasks consist in determining the auxiliary primal variables $(\bold{p},\bold{q})$, equations \ref{main4}b,c, and the auxiliary dual  variables $(\tilde{\bold{p}},\tilde{\bold{q}})$.
They are initialized to 0 and are updated as follows:
The primal variables $\bold{p}$ is updated through a TV proximity operator, which admits a closed form solution  given by equation \ref{primal_p} \citep[see][]{Combettes_2011_PRO} (step 3 of the algorithm). 
The primal variable $\bold{q}$ is updated by a projection operator, which also admits a closed form solution given by equation \ref{dual_q} (step 4 of the algorithm).
Finally, the duals are updated according to gradient ascend steps (step 5 of the algorithm)
\begin{equation}
\begin{cases}
\tilde{\bold{p}}^{k+1} &= \tilde{\bold{p}}^{k} + \bold{p}^{k+1}-\overline{\nabla} \bold{m}^{k+1}\\
\tilde{\bold{q}}^{k+1} &= \tilde{\bold{q}}^{k} + \bold{q}^{k+1}- \bold{m}^{k+1}.
\end{cases}
\end{equation}

%
\append{Using fourth-order equation for wavefield reconstruction ~ \small{(Step 1 of the algorithm)}} \label{Appb}
The wavefield reconstruction subproblem, equation \ref{UE}, can be written as the following over-determined system 
\begin{align} \label{4thU}
&\begin{bmatrix}
    \omega^2  \text{diag}(\bold{m}_{v_0})+ \text{diag}(\bold{m_\epsilon}) \nabla_{\!\! xx} & \text{diag}(\bold{m_\delta}) \nabla_{\!\! zz} \\
    \text{diag}(\bold{m_\delta}) \nabla_{\!\! xx} & \omega^2  \text{diag}(\bold{m}_{v_0}) + \nabla_{\!\! zz} \\
    \frac12\bold{\tilde{P}} & \frac12\bold{\tilde{P}}
\end{bmatrix}
\begin{bmatrix}
\bold{U}_x\\
\bold{U}_z\\
\end{bmatrix} \nonumber \\
&\hspace{4cm}=\begin{bmatrix}
\bold{S}_x\\
\bold{S}_z\\
\bold{D}
\end{bmatrix},
\end{align}
where $\bold{\tilde{P}}$ is the sampling operator of each component of the wavefield at receiver positions and $\bold{P}=[ \frac12 \bold{\tilde{P}}~~\frac12\bold{\tilde{P}}]$ (the coefficient $\frac12$ results because the (isotropic) pressure wavefield $\bold{u}$ recorded in the water is given by $\bold{u} = \frac12 \left( \bold{u}_x + \bold{u}_z\right)$). We can eliminate $\bold{U}_z$ from the equation \ref{4thU} to develop a fourth-order partial-differential equation for $\bold{U}_x$ and then update $\bold{U}_z$ from its explicit expression as a function of  $\bold{U}_x$ without any computational burden. \\
By multiplying the second row of equation \ref{4thU} by $\text{diag}(\bold{m}_\delta)$ and taking the difference with the first row, we find
\begin{align} \label{4thU1}
&\begin{bmatrix}
    \omega^2  \text{diag}(\bold{m}_{v_0})+ \text{diag}(\bold{m_\epsilon}) \nabla_{\!\! xx} & \text{diag}(\bold{m_\delta}) \nabla_{\!\! zz} \\
 \omega^2  \text{diag}(\bold{m}_{v_0})+ \text{diag}(\bold{m_\epsilon}-\bold{m}^2_{\delta}) \nabla_{\!\! xx} & -\omega^2\text{diag}(\bold{m_\delta} \circ \bold{m}_{v_0}) \\
    \frac12\bold{\tilde{P}} & \frac12\bold{\tilde{P}}
\end{bmatrix} \nonumber\\
&~~~~~~~~~~~~~~~~~~~
\begin{bmatrix}
\bold{U}_x\\
\bold{U}_z\\
\end{bmatrix}
=\begin{bmatrix}
\bold{S}_x\\
\bold{S}_x-\text{diag}(\bold{m_\delta})\bold{S}_z\\
\bold{D}
\end{bmatrix}.
\end{align}
where $\bold{m}^2_{\delta}=\bold{m_\delta} \circ \bold{m_\delta}$.
The second equation of \ref{4thU1} provides us the closed-form expression of $\bold{U}_z$ as a function of $\bold{U}_x$
\begin{equation}
\bold{U}_z=\bold{A}_z \bold{U}_x+\bold{B}_z,
\end{equation}
where
\begin{equation} \label{4thUz}
\begin{cases}
\bold{A}_z =  \frac{\text{diag}(\bold{m_\epsilon}-\bold{m}^2_{\delta}) \nabla_{\!\! xx}+ \omega^2\text{diag}(\bold{m}_{v_0})}
{\omega^2 \text{diag}(\bold{m}_{v_0} \circ \bold{m_\delta})} , \\
\bold{B}_z = \frac{\text{diag}(\bold{m_\delta})\bold{S}_z-\bold{S}_x}{\omega^2 \text{diag}(\bold{m}_{v_0} \circ \bold{m_\delta})}.
\end{cases}
\end{equation}
Injecting equation \ref{4thUz} into the first equation of \ref{4thU1} leads to the over-determined system satisfied by $\bold{U}_x$ 
\begin{align} \label{4thU5}
&\begin{bmatrix}
    \omega^2  \text{diag}(\bold{m}_{v_0})+ \text{diag}(\bold{m_\epsilon}) \nabla_{\!\! xx}+ \text{diag}(\bold{m_\delta}) \nabla_{\!\! zz}\bold{A}_z\\
    \frac12\bold{\tilde{P}} + \frac12\bold{\tilde{P}}\bold{A}_z
\end{bmatrix} \bold{U}_x \nonumber \\
&~~~~~~~~~~~~~~~~
=\begin{bmatrix}
\bold{S}_x-\text{diag}(\bold{m_\delta})\nabla_{\!\! zz}\bold{B}_z\\
\bold{D}-\frac12 \bold{\tilde{P}} \bold{B}_z
\end{bmatrix}.
\end{align}
So, instead of solving the $(n_r+2n) \times 2n$ linear system \ref{4thU} to update $\bold{U}$, $\bold{U}_x$ is updated by solving the $(n_r+n) \times n$ linear system \ref{4thU5} and then  $\bold{U}_z$ is updated using \ref{4thUz} without significant computational overhead. The wave equation operator in equation \ref{4thUz} has been broken down into an elliptic wave equation operator and anelliptic correction term (the term of $\bold{A}_z$ related to $\text{diag}(\bold{m}_{\epsilon}-\bold{m}_{\delta}^2)$, equation \ref{4thUz}). The former can be accurately discretized with the 9-point finite-difference method of \citet{Chen_2013_OFD}, while the anelliptic term can be discretized with a basic second-order accurate 5-point stencil without generating significant inaccuracies in the modelling \citep{Operto_2014_FAT}.
%
\newcommand{\SortNoop}[1]{}

\end{document}